\documentclass[10pt]{amsart}
\newtheorem{theorem}{Theorem}[section]
\newtheorem{lemma}[theorem]{Lemma}

\theoremstyle{definition}
\newtheorem{definition}[theorem]{Definition}
\newtheorem{example}[theorem]{Example}

\theoremstyle{remark}
\newtheorem{remark}[theorem]{Remark}

\numberwithin{equation}{section}

\theoremstyle{plain}
\newtheorem{corollary}[theorem]{Corollary}

\theoremstyle{remark}

\newtheorem{notation}[theorem]{Notation}

\usepackage{graphicx}
\usepackage{amssymb}

\begin{document}

\title{Valence of complex-valued
planar harmonic functions}

\author{Genevra Neumann}
\address{Department of Mathematics, 
University of California, Berkeley, California 94720}
\curraddr{Department of Mathematics, 
Kansas State University, 
Manhattan, Kansas 66506}
\email{neumann@math.ksu.edu}

\date{July 26, 2003;
revised January 26, 2004.}


\begin{abstract}
The valence of a function $f$ at a point $w$
is the number of distinct, finite solutions to
$f(z) = w$.
Let $f$ be a complex-valued harmonic function in an 
open set $R \subseteq \mathbb{C}$.
Let
$S$ denote the critical set of $f$
and $C(f)$ the global cluster set of $f$.
We show that
$f(S) \cup C(f)$ partitions
the complex plane into regions of constant valence.
We give some conditions such that $f(S) \cup C(f)$
has empty interior.
We also show that a component 
$R_0 \subseteq 
R \backslash f^{-1} (f(S) \cup C(f))$ is
a $n_0$-fold covering of some component
$\Omega_0 \subseteq
\mathbb{C} \backslash (f(S) \cup C(f))$.
If $\Omega_0$ is simply connected,
then $f$ is univalent on $R_0$.
We explore conditions for
combining adjacent components
to form a larger region of univalence.
Those results which hold for $C^1$ functions
on open sets in $\mathbb{R}^2$ are first stated
in that form and then applied to the case of
planar harmonic functions.
If $f$ is a light, harmonic function in the
complex plane,
we apply a structure theorem of Lyzzaik
to gain information about the difference in valence
between components of 
$\mathbb{C} \backslash (f(S) \cup C(f))$
sharing a common boundary arc in 
$f(S) \backslash C(f)$.
\end{abstract}

\maketitle

\section{Introduction}
\label{sec:intro}
We study here complex-valued harmonic functions in 
the plane,
which we refer to simply as harmonic functions.
The behavior of such a function can be vastly 
different from that of a holomorphic function.
Analytic polynomials take every value a finite number
of times.
In contrast, the range of a harmonic polynomial
can exclude an open region of the complex plane.
Picard's theorem states
that a transcendental entire function takes every value,
with the exception of possibly one point, an infinite
number of times.
In contrast, there are transcendental harmonic 
functions that omit open regions.
Also, there are transcendental harmonic functions
that approach $\infty$ as $z \rightarrow \infty$
(like analytic polynomials),  
such that each $w \in \mathbb{C}$ has
a finite number of distinct preimages
(like analytic polynomials), 
and such that the maximum possible number of
preimages is unbounded 
(like transcendental entire functions.)

A harmonic polynomial $f(z)$ is a harmonic function
of the form
$f(z) = p(z) + \overline{q(z)}$
where $p$ and $q$ are analytic polynomials in $z$.
Let $n_p$ be the degree of $p$ as a polynomial in $z$
and $n_q$ be the degree of $q$.
A. Wilmshurst~\cite{W:thesis, W:paper}\footnote{We 
will refer to~\cite{W:paper} 
when a result appears in both~\cite{W:thesis}
and~\cite{W:paper}.}
showed that when 
$f$ has a finite number of zeros,
it has at most $N^2$ distinct zeros, 
where $N = max(n_p, n_q)$.
Wilmshurst's bound is sharp;
there are examples~\cite{BHS:maxval, W:paper} of 
harmonic polynomials with 
$N^2$ distinct zeros when $n_p = N$ and $n_q = N-1$.
For $1 \leq n_q \leq n_p - 1$,
Wilmshurst conjectured that $f$ has 
at most $n_q (n_q - 1) + 3n_p - 2$ distinct zeros.
D. Khavinson and G. \'{S}wi\c{a}tek~\cite{KS:3n-2} 
recently proved Wilmshurst's conjecture for the 
case $n_q = 1$
using methods from complex dynamics.

We will be concerned with questions related to
the valence of harmonic functions.
The valence of a function $f$ at a given point $w$,
denoted $Val(f,w)$,
is the number of distinct points $z$ in the domain
of $f$ such that $f(z)=w$.  
The valence of a function, denoted $Val(f)$,
is the supremum of $Val(f,w)$.

Let $f(z) = (u(z), v(z))$ be a $C^1$ function in
an open set
$R \subseteq \mathbb{R}^2$.
The Jacobian of $f$ is given by 
$J_f = u_x v_y - u_y v_x$.
The inverse function theorem tells us that
if $J_f(z) \neq 0$,
then $f$ is a homeomorphism in 
some neighborhood of $z$.
H. Lewy~\cite{Lew:jacobian} proved the converse
when $u$ and $v$ are both real-valued harmonic
functions.
Hence,
if $f = u + iv$ is harmonic in an open set 
$R \subseteq \mathbb{C}$,
then $f$ is a local homeomorphism at $z$ if and only if
$J_f(z) \neq 0$.
The critical set $S$ of a $C^1$ function $f$ 
consists of 
those points where the Jacobian vanishes.
If $f$ is harmonic,
the critical set consists of those points where
$f$ is not locally 1-1.

When studying specific examples of harmonic
polynomials in $\mathbb{C}$ and looking at graphs
of the image of the critical set using 
Mathematica\footnote{Mathematica is
a registered trademark of Wolfram Research, Inc.},
the author noticed that the image of the critical set
partitions the complex plane into regions
of constant valence.
This is not true for all harmonic functions.
If $\lim_{z \rightarrow \infty} |f(z)| \neq \infty$,
then we need to include another set in order 
to partition
the complex plane into regions of constant valence.

Let $C(f, \infty)$ denote those finite values $w$
such that we can find a sequence $\{z_n\}$ with
$z_n \rightarrow \infty$ and 
$f(z_n) \rightarrow w$;
$C(f, \infty)$ is the 
cluster set of $f$ at $\infty$.
The global cluster set of a function $f$ 
with domain $R$ is denoted
$C(f)$ and consists of all 
finite values which are approached on some sequence 
of points in $R$ which converge to a point in 
$\partial R \cup \{\infty\}$.
(Our definitions are slightly different from those in 
the book of E. Collingwood and 
A. Lohwater~\cite{CL:cluster},
where a cluster set can include the point
at infinity.
We are interested in partitions of the finite
complex plane and exclude the point at infinity.)

We show, 
under suitable conditions,
that the plane can be partitioned
into regions of constant valence.
Since many of our results hold not only for 
harmonic functions but for $C^1$
functions in open subsets of $\mathbb{R}^2$,
we will state the partitioning results first for
$C^1$ functions defined in an open set in 
$\mathbb{R}^2$ and then specialize 
to the case where the function is harmonic
in an open set in $\mathbb{C}$.
We will then look at the case when $f$ is harmonic
in all of $\mathbb{C}$.
With each additional assumption on $f$,
we can say more about the partitioning set.
When $f$ is a light harmonic function,
we can apply some results of A. Lyzzaik~\cite{L:light}
to show that the valences for two regions
separated by an arc in the image of the critical set
differ by a non-zero, even number,
if that arc contains some point not in the cluster
set.

Let $R_1$ be a component of 
$R \backslash (f^{-1}(f(S) \cup C(f)))$.
We will see that $f(R_1)$ is a component of
$\mathbb{R}^2 \backslash (f(S) \cup C(f))$ and
that $f|_{R_1}$ is an even cover of $f(R_1)$
in the sense of Munkres~\cite[p. 331]{Mun:top}.
Suppose that $R_2$ is another component of this
partition of $R$;
suppose also that $R_2$ shares a common boundary arc
with $R_1$.
When $f$ is harmonic,
we will look at conditions for $R_1$ and $R_2$
to be mapped to different
components of $\mathbb{C} \backslash (f(S) \cup C(f))$.
We will also explore the behavior of $f$ on 
$\partial R_1$ 
(including the behavior at puncture points.)
If $f$ is univalent in $R_1$ and in $R_2$,
can we join $R_1$ and $R_2$ along the interior of their
shared boundary arc to get a larger region
of univalence? 

The figures below were produced from EPS files
generated by
Mathematica routines\footnote{These routines are
enhanced versions of the Mathematica routines 
in \cite{gcn:thesis}.
For the figures in this paper,
$z \in f^{-1}(w)$ if $|Re\,f(z) - Re\,w| + 
|Im\,f(z) - Im\,w|
< 10^{-10}$.} 
written by the author.
Some of the figures contain alphanumeric labels;
these labels were manually inserted into the
EPS files.

\section*{Acknowledgements}
Many of the results in this paper appeared in
the author's thesis~\cite{gcn:thesis}.
The author is grateful to her advisor,
D. Sarason,
for encouraging her to study planar harmonic functions.
She also wishes to thank him for his open-mindedness,
patience, 
and kindness in helping
her learn about mathematics and research.
His many suggestions were essential to this work.
The author is indebted to W. Hengartner for
suggesting that her earlier results concerning
regions of constant valence should be
extendable to the case of infinite valence,
and to H. Helson for suggesting that
they should be extendable to the case of functions
in arbitrary open regions
(not just in the entire plane.)
The author is also indebted to D. Khavinson for
information concerning the Brelot-Choquet Lemma.
The author also wishes to thank A. Bakke for
adding labels to some of the EPS files generated by 
Mathematica.

\section{Notation and background material}
\label{sec:notation}
Let $R$ be an open set in $\mathbb{R}^2$.
Let $f: R \rightarrow \mathbb{R}^2$ be $C^1$.
We can also write $f$ as 
$f(z) =(u(z),\ v(z))$, where $u,\ v:
R \rightarrow \mathbb{R}$
and $z=(x,y) \in \mathbb{R}^2$.
Then

\[
\begin{array}{rll}
J_f(z) &= &u_x(z)v_y(z)-u_y(z)v_x(z)\\
S &= &\{ z \in R: J_f(z) = 0 \}\\
C(f) &= &\{ \zeta \in \mathbb{R}^2 :
\exists \{ z_n \}_{n=1}^{\infty} \subset R
\textnormal{\ with\ }
\lim_{n \rightarrow \infty} z_n 
\in \partial R \cup \{\infty\}\ $and$ \\ 
&&\lim_{n \rightarrow \infty} f(z_n) = \zeta\}\\
C(f, \infty) &= &\{ \zeta \in \mathbb{R}^2 :
\exists \{ z_n \}_{n=1}^{\infty} \subset \mathbb{R}^2
\textnormal{\ with\ }
\lim_{n \rightarrow \infty} |z_n| \rightarrow 
\infty\ $and$ \\
&&\lim_{n \rightarrow \infty} f(z_n) = \zeta\}\\
B(w, \epsilon) &= &\{z:\  |z-w| \  < \epsilon\}\\
Val(f, w) &= &\#\{z \in R: f(z) = w\} \\
V_N(f) &= &\{w: w \notin f(S) \cup C(f) 
\textnormal{\ and\ } Val(f, w)=N\}\\
Val(f, V) &= &sup_{w \in V} \{Val(f, w)\} 
\end{array}
\]

\begin{enumerate}
\item $S$ denotes the critical set of $f$.
\item $C(f)$ denotes the global cluster set in the 
finite plane of $f|_R$.
Notice that if $R = \mathbb{R}^2$,
then $C(f) = C(f, \infty)$.
\item
$ \overline{f(S)} \subseteq f(S) \cup C(f) $:
Let $\{w_n\}_{n = 1}^ \infty \subseteq f(S)$.
Suppose that $\lim_{n \rightarrow \infty} w_n = w_0$,
where $w_0$ is finite.
For each $w_n$,
there exists $z_n \in S$ such that $f(z_n) = w_n$.
If $\{z_n\}_{n = 1} ^ \infty$ is unbounded,
then $w_0 \in C(f)$ by definition.
If $\{z_n\}_{n = 1} ^ \infty$ is bounded,
then it has a convergent subsequence,
say $\{z_m\}$,
such that $z_m \rightarrow z_0 \in \mathbb{R}^2$.
If $z_0 \in R$,
then $z_0 \in S$ by the continuity of $J_f(z)$;
hence $w_0 \in f(S)$.
Otherwise,
we must have that $z_0 \in \partial R$;
hence $w_0 \in C(f)$.

\item
It is well known that $C(f)$ is closed.

\item
$ f(S) \cup C(f) $ is closed,
since $\overline{f(S)} \subseteq f(S) \cup C(f)$.

\item
$ f^{-1} \left( f(S) \cup C(f) \right) $ 
is relatively closed in $R$, 
since $f$ is continuous.
\end{enumerate}

\subsection {Lyzzaik's local description of 
light harmonic functions}
\label{sec:lyzintro}
Recall that a function is said to be light if
the preimage of each point is empty or totally
disconnected.
A. Lyzzaik~\cite{L:light} has characterized the 
local behavior of a function $f$ which is light and
harmonic in a simply connected, open set in
$\mathbb{C}$.
Lyzzaik pays special attention to the behavior
of $f$ on the critical set and the behavior of 
$f$ in a neighborhood of the critical set.
Let $S$ denote the critical set of $f$.

We first note that since $f$ is harmonic in a 
simply connected, open set,
we can find functions $h$ and $g$ holomorphic
in that set such that $f(z) = h(z) + \overline{g(z)}$.
Thus,
$z \in S$ if and only if 
$0 = J_f(z) = |h'(z)|^2 - |g'(z)|^2$.
Lyzzaik defines the meromorphic function 
$\psi = h' / g'$ and uses this function to
study the local behavior of $f$ near the critical set.
Note that if $z \in S$, 
then either
$\psi(z)$ is unimodular,
$\psi(z) = 0$,
or $\psi$ has a pole at $z$.
Since $f$ is also light,
$S$ has empty interior and 
$h'(z) = 0 = g'(z)$ at isolated points
(see Lemma \ref{lem:nullpartifsnotempty}.)
Hence the zeros and poles of $\psi$ 
which lie in $S$ are isolated.
In Lyzzaik's classification of critical points,
the set of isolated critical points is denoted $N$.
Since $\psi$ is unimodular at the remaining points
in $S$,
if $z \in S \backslash N$,
then $z$ lies in an analytic arc in $S$.
The set of branch points of $S$ is denoted $F_3$.
Hence,
if $z \in S \backslash (N \cup F_3)$,
in some small neighborhood of $z$,
there is a unique analytic arc in $S$ with $z$
in its interior.
Lyzzaik also classifies those critical points
$z \in S \backslash (N \cup F_3)$ where
the the image of $f$ stops.
Let $\gamma \subset S \backslash (N \cup F_3)$ 
be an analytic
arc in $S$ such that $z \in int\ \gamma$.
Suppose that $f(\gamma)$ has zero speed at $f(z)$.
If the argument of the tangent to $f(\gamma)$
jumps by $\pm \pi$ at $f(z)$,
then $f$ is said to have a harmonic cusp at
$f(z)$ and $z \in F_1$.
Otherwise,
if $z \in S \backslash (N \cup F_1 \cup F_3)$ and
$f'(z) = 0 = g'(z)$, 
then $z \in F_2$.

Lyzzaik shows that 
$N \cup (\cup_{j=1}^3 F_j)$ consists of isolated 
points,
and if 
$\gamma \subset S \backslash (N \cup F_3)$,
then $f|_\gamma$ is a local homeomorphism.
Further, 
if $z \in S \backslash (N \cup F_1 \cup F_3)$,
then a subarc of $S$ with $z$ in its interior
is mapped to a convex arc\footnote{Lyzzaik defines
a convex arc to be 
a directed simple arc where
the slope of the tangent is 
continuously increasing.} 
by $f$.
Lyzzaik also shows that if $z \in F_3$,
then $f(S)$ in some neighborhood of $z$ consists
of convex arcs and harmonic cusps.
Lyzzaik's classification of critical points
is discussed in more detail in 
Section \ref{sec:lyzzaik} below.

We will also use Lyzzaik's 
local characterization
of $f$ in a neighborhood of $z_0$ for
$z_0 \in S \backslash (N \cup F_3)$.
This result is given below as Theorem
\ref{thm:lyzthm5.1}.
We now review some background for
this result and explain the notation.

Let $D$ denote the unit disc centered
at the origin.
A function $f$ is said to be locally topologically
$z^n$ at $z_0$,
denoted $f_{z_0} \sim z^n$,
if there exist an open neighborhood $U$ of $z_0$
and homeomorphisms $h_1: U \rightarrow D$
and $h_2: \mathbb{C} \rightarrow \mathbb{C}$ such that
$h_2 \circ f(z_0) = 0$ and such that
$h_2 \circ f \circ h_1^{-1}(\zeta) = \zeta^n$
for all $\zeta \in D$.
A result of S. St\"{o}ilow 
(see~\cite{Sto:lecons} or Lemma 1 in~\cite{AL:local})
says that if $f$ is continuous and open in a
neighborhood of $z_0$,
then $f$ is locally topologically $z^n$ for some
positive integer $n$.
Y. Abu-Muhanna and A. Lyzzaik~\cite{AL:local} 
extend St\"{o}ilow's result to points in the
boundary:

\begin{lemma}[Abu-Muhanna and Lyzzaik]
\label{lem:ALlem2}
Let $G^+ = \{z: |z| < 1, Im\ z > 0\}$ 
be the open semi-disc and
$G = \{z: |z| < 1, Im\ z \ge 0\}$ be
the half-closed semi-disc.
Suppose that $f: G \rightarrow \mathbb{C}$ is a 
continuous function open in $G^+$ and 
topological
on $G \backslash G^+$.
Then for every $z_0 \in G \backslash G^+$
there is a positive integer $n$ such that 
$f$ at $z_0$ is locally topologically $z^{2 n - 1}$. 
This result also holds for sets homeomorphic 
to $G^+$ and $G$.
\end{lemma}

\noindent
They prove this by constructing a function $F(z)$
using an idea similar to that used in proving the
Schwarz reflection principle and showing that
$F$ is an open map on a disc centered at $z_0$.
Then they apply St\"{o}ilow's result to $F$.
The result follows by restricting $h_1$ and
$F$ to $G$.

Suppose that $f$ is a harmonic function in some
open set $R \subseteq \mathbb{C}$.
Given $z_0 \in R \backslash S$,
$f$ is either sense-preserving ($J_f(z) > 0$)
or sense-reversing ($J_f(z) < 0$) in some open
neighborhood of $z_0$.
Suppose that $f_{z_0} \sim z^n$.
If we require $h_1$ and $h_2$ to be sense-preserving
homeomorphisms,
we see that $f_{z_0} \sim z^n$ if $f$ is 
sense-preserving in a punctured neighborhood $z_0$
and that $f_{z_0} \sim \overline{z}^n$ if 
$f$ is sense-reversing in a punctured neighborhood
of $z_0$.

We now suppose that $f$ is a light harmonic function
in some simply connected, open set.
Suppose that $z_0 \in int\ \gamma_0 \subset
S \backslash (F_3 \cup N)$.
Then $z_0$ can be thought of as living in a shared
boundary arc of two adjacent regions.
In a sufficiently small neighborhood $U$ of
$z_0$,
$f$ is 1-1 on $\gamma = \gamma_0 \cap U$ 
(recall that $f|_\gamma$ is a local homeomorphism at
$z_0$) and 
$\gamma$
splits $U$ into a sense-preserving region $R^+$
and a sense-reversing region $R^-$.
Lyzzaik~\cite{L:light} notes that 
Lemma \ref{lem:ALlem2}
applies to $R^+ \cup \gamma$
and to $R^- \cup \gamma$.
The notation $f_{z_0} \sim z^j, \overline{z}^k$
means that  
$f$ at $z_0$
is locally topologically $z^j$ in $R^+ \cup \gamma$
and that $f$ at $z_0$
is locally topologically $\overline{z}^k$ 
in $R^- \cup \gamma$.
Lemma \ref{lem:ALlem2} also shows that $j$ and $k$ must
both be odd positive integers.
A structure theorem of Lyzzaik for this case
(Theorem \ref{thm:lyzthm5.1} below) 
gives the values for $j$ and $k$
based on various conditions on $z_0$.

In particular,
suppose that $z_0 \in S \backslash (N \cup 
(\cup_{i=1}^3 F_i))$.
Then Theorem \ref{thm:lyzthm5.1}
gives $f_{z_0} \sim z, \overline{z}$
since $g'(z_0) \neq 0$.
We can assign a direction to $f(\gamma)$,
which is part of the common boundary of $f(R^+)$ and
$f(R^-)$.
We see that the image of $U$ is folded over
$f(\gamma)$;
in other words,
$f(U \backslash \gamma)$ lies to one side of 
$f(\gamma)$.
The tangent line to $f(\gamma)$
at $f(z_0)$ lies in $f(U)$.
One way to see this is to recall that 
$f$ is univalent and sense-preserving in $R^+$ 
($f_{z_0} \sim z$ for $z \in R^+$.)
Also,
$f$ is locally 1-1 on $\gamma$.
A result of P. Duren and D. 
Khavinson~\cite{DK:concave} shows that
$f(\gamma)$ is concave with respect to $f(R^+)$;
hence the tangent line to $f(\gamma)$ at $f(z_0)$
lies in $f(U)$.

In Remark \ref{rem:sharedarc},
we note that Lyzzaik's result implies that
two components of $\mathbb{C} \backslash
f^{-1} (f(S) \cup C(f, \infty))$ 
which share a common boundary
arc in $S$ must be mapped to the same component of
$\mathbb{C} \backslash (f(S) \cup C(f, \infty))$.
In Section \ref{sec:lyzzaik},
we use Lyzzaik's structure theorem to compare
the valence in two components of
$\mathbb{C} \backslash (f(S) \cup C(f, \infty))$ when
the two components share a common boundary arc
in $f(S)$.

\section{Partitioning results for $C^1$ mappings
in $\mathbb{R}^2$}
\label{sec:c1part}
We show that $f(S) \cup C(f)$ partitions
$\mathbb{R}^2$ into regions of constant valence
when $f: R \rightarrow \mathbb{R}^2$ is $C^1$ 
in an open set $R \subseteq \mathbb{R}^2$.
We then examine a corresponding partition of
$R$.

\subsection {Regions of constant valence for 
$\mathbf{C^1}$ mappings in $\mathbb{R}^2$}

\begin{lemma}
\label{lem:c1lsc}
Let $R$ be open in $\mathbb{R}^2$.
Let $f: R \rightarrow \mathbb{R}^2$ be
$C^1$.
Let $w_0 \in \mathbb{R}^2 \backslash \overline{f(S)}$.
Suppose that $Val(f, w_0) \ge N_0 \ge 0$,
where $N_0$ is finite.
Then there exists an open neighborhood of $w_0$,
say $W_0$,
such that $Val(f, w) \ge N_0$ for all $w \in W_0$.
\end{lemma}
\begin{proof}
Trivial if $\overline{f(S)} = \mathbb{R}^2$, 
so suppose that $\mathbb{R}^2 \backslash 
\overline{f(S)} \neq \varnothing$.
The result is also trivial if $N_0 = 0$,
so we will suppose that $N_0 > 0$.
Since $Val(f, w_0) \ge N_0$,
choose $N_0$ distinct points $z_1, ..., z_{N_0}$
in $f^{-1}(w_0)$.
By assumption, these points will be in the open set 
$R \backslash f^{-1}(\overline{f(S)})$.
By the inverse function theorem,
we may find an open neighborhood of each $z_j$, 
say $U_j$,
such that $f: U_j \rightarrow W_j$ is 1-1 and onto,
where $W_j$ is an open neighborhood of $w_0$.
Since $N_0$ is finite,
we can choose the $U_j$ to be pairwise disjoint.
Let $W_0 = \cap_{j=1}^{N_0} W_j$.
Then $W_0$ is an open, non-empty set such that 
each $w \in W_0$
has a distinct preimage in each $U_j$.
The lemma then follows.
\end{proof}

\begin{lemma}
\label{lem:nvalopen}
Let $R$ be an open set in $\mathbb{R}^2$.
Let $f: R \rightarrow \mathbb{R}^2$ be $C^1$.
Fix $N_0 \ge 0$, finite. 
Then $V_{N_0}(f)$ is open in $\mathbb{R}^2$.
\end{lemma}
\begin{proof}
Let $ Q_0 = \mathbb{R}^2 \backslash (f(S) \cup 
C(f)) $.
Then $ Q_0 $ is open.
By definition, $V_{N_0}(f) \subseteq Q_0$.
If $ V_{N_0}(f) = \varnothing $, the claim is vacuously 
true.
Suppose that $ V_{N_0}(f) \neq \varnothing $;
hence we also have $Q_0 \neq \varnothing$.
Choose $w_0 \in V_{N_0}(f)$.
Thus, $Val(f,w_0) = N_0$ and $w_0 \in Q_0$.

We first note that there exists an 
open neighborhood $ \tilde{B} \subseteq Q_0$
of $w_0$ such that $Val(f,w) \ge N_0$ for 
all points
$w \in \tilde{B}$.
Since $\overline{f(S)} \subseteq f(S) \cup C(f)$,
this follows from Lemma \ref{lem:c1lsc}.

Since $ \tilde{B}$ is open, 
we may choose $\epsilon_0 > 0$ such that 
$B(w_0, \epsilon_0) \subseteq \tilde{B}$.
Suppose that we can find $\{\epsilon_j\}_{j=1}^{\infty}$ 
with
$\lim_{j \rightarrow \infty}
\epsilon_j = 0$,
such that for each $j$, $0 < \epsilon_j < \epsilon_0$
and there exists 
$w_j \in B(w_0, \epsilon_j)$ with $Val(f, w_j) > N_0$.
By construction, $w_j \rightarrow w_0$ and
$Val(f, w_j) \geq N_1$ where $N_1 = N_0 + 1$.
Without loss of generality, 
we may suppose that the $w_j$ are distinct.
We will show that $Val(f,w_0) \geq N_1 > N_0 = 
Val(f,w_0)$,
a contradiction.

We may find
$ z_{j1}, z_{j2}, \ldots, z_{jN_1} \in R 
\backslash f^{-1} \left( f(S) \cup C(f) \right) $, 
distinct, such that $f(z_{jk}) \\
= w_j$ 
for $ k = 1, 2, \ldots, N_1 $,
since $w_j \notin f(S) \cup C(f)$ and
$Val(f, w_j) \geq N_1$.
Moreover,
$ \bigcup_{k=1}^{N_1} 
\{ z_{jk} \}_{j=1}^{\infty} $ 
consists of pairwise distinct points,
since the $ w_j $ are distinct.

\begin{enumerate}

\item
$ \bigcup_{k=1}^{N_1} 
\{ z_{jk} \}_{j=1}^{\infty} $ is bounded: 
This is obvious if $R$ is bounded.
So, suppose that $R$ is unbounded. 
Suppose that $\bigcup_{k=1}^{N_1} 
\{ z_{jk} \}_{j=1}^{\infty} $ is unbounded.
Thus, we may find $ \{ z_{j_lk_l} \}_{l=1}^{\infty} $ 
where $ k_l \in \{ 1, 2, \ldots, N_1 \} $
such that $ |z_{j_lk_l}| \rightarrow \infty $ as 
$ l \rightarrow \infty $.
By the pigeonhole principle, we may find a 
subsequence where $ k_l $ is fixed; \textit{i.e.}, 
$ \{ z_{lN} \}_{l=1}^{\infty} $
for some \mbox{$ N \in \{ 1, 2, \ldots, N_1 \} $}
such that $ |z_{lN}| \rightarrow \infty $ as 
$ l \rightarrow \infty $.
Since $ \{ z_{lN} \}_{l=1}^{\infty} 
\subseteq \{ z_{jN} \}_{j=1}^{\infty} $ and
$ \lim_{j \rightarrow \infty} 
f(z_{jN}) = \lim_{j \rightarrow \infty} 
w_j = w_0 $,
$ \lim_{l \rightarrow \infty} 
f(z_{lN}) = w_0 $.
Since $ |z_{lN}| \rightarrow \infty $, $ w_0 \in 
C(f) $.
But $ C(f) \cap Q_0 = \varnothing $ and $ w_0 
\in Q_0 $ by assumption.
Thus $ \bigcup_{k=1}^{N_1} 
\{ z_{jk} \}_{j=1}^{\infty} $ is bounded.

\item
If $z$ is a finite cluster point of 
$ \bigcup_{k=1}^{N_1} 
\{ z_{jk} \}_{j=1}^{\infty}  $, then $z \in R$
and $ f(z) = w_0 $:
Find $ \{ z_{j_lk_l} \}_{l=1}^{\infty} $ 
where $ k_l \in \{ 1, 2, \ldots, N_1 \} $ such that
$ \lim_{l \rightarrow \infty} 
z_{j_lk_l} = z $.
As above, we may find a subsequence of this sequence 
with $ k_l $ constant, say 
$ \{ z_{lN} \}_{l=1}^{\infty}$ 
such that $z_{lN} \rightarrow z$.
Since $ \{ z_{lN} \}_{l=1}^{\infty} 
\subseteq \{ z_{jN} \}_{j=1}^{\infty} $
and $ f(z_{jN}) \rightarrow w_0 $, 
$ \lim_{l \rightarrow \infty} 
f(z_{lN}) = w_0 $.
Since $z$ is a cluster point of a subset of $R$,
either $z \in R$ or $z \in \partial R$.
If $z \in \partial R$,
we have $w_0 \in C(f)$,
a contradiction.
Hence $z \in R$.
Since $ z_{lN} \rightarrow z $ and $f$ is continuous
at $z$, 
$ f(z) = w_0 $.

\item
$ \bigcup_{k=1}^{N_1} 
\{ z_{jk} \}_{j=1}^{\infty}  $ has at least 
$ N_1 $ distinct cluster points in $R$:
By (1), all of the cluster points are finite.
By (2),
all of these cluster points are in $R$.
Let $z$ be a cluster point.
Suppose that for each $ \epsilon > 0 $, there is a 
$ j>0 $ such that $ |z_{jk_1} - z| < \epsilon $ 
and $ |z_{jk_2} - z| < \epsilon $
for some choice of $ k_1, k_2 $ with $ k_1 \neq k_2 $.
But for each such $j$, $ f(z_{jk_1}) = w_j = 
f(z_{jk_2}) $, with $ z_{jk_1} \neq z_{jk_2} $,
so $f$ is not locally 1-1 at $z$; hence $ z \in S $.
By (2), $ f(z) = w_0 \in Q_0 $ with $ Q_0 
\cap f(S) = \varnothing $, so $ z \notin S $, 
a contradiction.
So, there exists $ \epsilon > 0 $ where for each 
$ j > 0 $ such that $ | z_{jk} - z | < \epsilon $ 
holds, it holds for exactly one value of $k$
for that choice of $j$.
Thus we must have at least $ N_1 $ cluster points
in $R$.

\end{enumerate}

\noindent
From (1)-(3), we see that $ w_0 $ has at least 
$ N_1 $ distinct preimages in $R$,
which gives the desired contradiction.
With this contradiction, 
we have shown that $\exists\ \epsilon > 0$ such
that $Val(f, w) = N_0$ for all $w \in B(w_0, \epsilon)$.
Since $w_0 \in V_{N_0}(f)$ is arbitrary,
$V_{N_0}(f)$ is open in 
$\mathbb{R}^2$. 
\end{proof}

\begin{lemma}
\label{lem:infval}
Let $f$ be a $C^1$ mapping defined in an open set 
$R \subseteq \mathbb{R}^2$.
If $Val(f, w_0) = \infty$,
then $w_0 \in f(S) \cup C(f)$.
\end{lemma}
\begin{proof}
Since $Val(f, w_0) = \infty$,
we may choose a sequence $\{z_n\} \subseteq 
f^{-1}(w_0)$ consisting of distinct points.
If $\{z_n\} \subseteq f^{-1}(w_0)$
has a bounded subsequence 
(which we also denote by $\{z_n\}$)
converging to a point
$z^* \in R$,
we will show that $z^* \in S$.
Given $\epsilon > 0$,
we can find $N$ such that $|z_n - z^*| < \epsilon$
for all $n > N$.
By continuity,
$f(z^*) = w_0$.
Hence $f$ is not locally 1-1 at $z^*$.
By the inverse function theorem,
$J_f(z^*) = 0$.
Hence $z^* \in S$ and $w_0 = f(z^*) \in f(S)$.

Otherwise,
we may suppose that either $\{z_n\}$ has a subsequence
converging to a finite point in $\partial R$ or
has an unbounded subsequence.
In either case,
since each point of the subsequence is mapped to
$w_0$,
we have $w_0 \in C(f)$.
\end{proof}
\textbf{Comment:}
It would be nice to know if the following result is similar to known
results in differential topology or of other applications of similar
results.
\begin{theorem}
\label{thm:c1part}
Let $R$ be open in $\mathbb{R}^2$. 
Let $f: R \rightarrow \mathbb{R}^2$ be $C^1$.
Then $f(S) \cup C(f)$ partitions 
$\mathbb{R}^2$ into regions of constant valence.
\end{theorem}
\begin{proof}
Let $\varphi(w) = Val(f, w)$.
It is enough to show that $\varphi$ is a continuous,
integer-valued function at each $w \in
\mathbb{R}^2 \backslash (f(S) \cup C(f))$.
Choose $w_0 \in
\mathbb{R}^2 \backslash (f(S) \cup C(f))$.
By Lemma \ref{lem:infval}, 
$\varphi(w_0)$ is a finite integer,
say $j$.
Hence $w_0 \in V_j$.
By Lemma \ref{lem:nvalopen},
$V_j$ is open.
By definition,
$V_j \subseteq \mathbb{R}^2 \backslash
(f(S) \cup C(f))$.
Hence there exists $\delta > 0$ such that
$B(w_0, \delta) \subseteq V_j$ and
$\varphi(w) = j$ for all $w \in B(w_0, \delta)$.
Thus $\varphi$ is a continuous, integer-valued
function in every region off $f(S) \cup C(f)$
and the conclusion follows.  
\end{proof}

From Lemma \ref{lem:c1lsc} and
Theorem \ref{thm:c1part},
we see that $Val(f,w)$ is lower semi-continuous
on $\mathbb{R}^2 \backslash \overline{f(S)}$
for points with finite valence.

\begin{example}
\label{ex:c1poly}
$f(x,y) = (x^2 + y^2, 2xy)$

Here $R = \mathbb{R}^2$,
so $C(f) = C(f, \infty)$.
Clearly $|f| \rightarrow \infty$ as 
$z \rightarrow \infty$,
so $C(f, \infty) = \varnothing$.
The critical set consists of the lines
$y = x$ and $y = -x$.
$f$ maps the critical set to the rays $y = x$ 
and $y = -x$ for $x \ge 0$.
A calculation shows that $f(S)$ partitions
$\mathbb{R}^2$ into regions of constant valence.
In particular, each point with $x > 0$
and $|y| < x$ has four distinct preimages. 
The origin has one preimage.
Each point in the image of the critical set 
in the right half plane 
has two preimages.
The remaining points have no preimages.
Note that 
the behavior of $Val(f,w)$ on the partitioning
set is consistent 
with $Val(f,w)$ being lower semi-continuous on 
$\mathbb{R}^2 \backslash \overline{f(S)}$.
If we rewrite $f$ as $f = u + iv$ where $z = x + iy$, 
it is clear that $f(z)$ is not harmonic.
\end{example}

We can say a little about the behavior of $Val(f,w)$
when $Val(f,w)$ is infinite.

\begin{lemma}
\label{lem:c1infvallsc}
Let $R$ be open in $\mathbb{R}^2$.
Let $f: R \rightarrow \mathbb{R}^2$ be
$C^1$.
Suppose that $Val(f, w_0) = \infty$ and that
$f^{-1}(w_0) \backslash S$ contains an infinite
number of distinct points.
Then, given $N_0 \ge 0$, finite, 
there exists an open neighborhood of $w_0$,
say $W_0$,
such that $Val(f, w) \ge N_0$ for all $w \in W_0$.
\end{lemma}
\begin{proof}
Obvious if $N_0 = 0$, so assume $N_0 > 0$.
Choose $\{z_1, ..., z_{N_0}\} \subseteq
f^{-1}(w_0) \backslash S$, distinct.
By the inverse function theorem,
for each $z_j$, there exists an open neighborhood
of $z_j$, say $B_j$, such that $f$ is 1-1 on $B_j$
and $f(B_j)$ is open.
We may choose the $B_j$ so that 
they are pairwise disjoint.
Let $W_0 = \cap f(B_j)$.
\end{proof}

Example \ref{ex:flatpoly} below demonstrates
why we must require that $w_0$ have an infinite
number of distinct preimages off of the critical set.
In this example,
the origin has infinite valence and every neighborhood
of the origin contains a point with no preimages.
However,
the preimages of the origin all lie in the critical 
set.

\begin{theorem}
\label{thm:harmlsc}
Let $f$ be a $C^1$ mapping in an open set
$R \subseteq \mathbb{R}^2$.
Let $V$ be a connected component of 
$\mathbb{R}^2 \backslash (f(S) \cup C(f))$
and let $N_0 = Val(f, V)$.
Let $w_0 \in \partial V \backslash \overline{f(S)}$.
Then $Val(f, w_0) \le N_0$.
\end{theorem}
\begin{proof}
By Theorem \ref{thm:c1part} and Lemma \ref{lem:infval},
$f$ has constant, finite valence on $V$;
hence $N_0$ is finite.
Suppose that $Val(f, w_0) > N_0$.
Then,
by Lemma \ref{lem:c1lsc}, 
there exists an open neighborhood $V_0$
of $w_0$ such that $Val(f, w) > N_0$ for all
$w \in V_0$.
Since $w_0 \in \partial V$,
$V_0 \cap V$ is a non-empty open set.
For $w \in V_0 \cap V \subseteq V$,
$Val(f, w) > N_0 = Val(f, w)$, 
a contradiction and the theorem follows.
\end{proof}

The preceding two results show that $Val(f,w)$ is
lower semi-continuous off of $\overline{f(S)}$.
Consider $f(z) = z + Re\ e^z$ in 
Example \ref{ex:transharm} below.
Since $R = \mathbb{C}$,
$C(f) = C(f, \infty)$.
Each $w \in f(S) \cup C(f, \infty)$
has exactly one preimage.
Each point in $f(S)$ has a neighborhood containing
points with no preimages and points with two
preimages,
so this example shows why $f(S)$ is excluded
in the preceding result.
If $w_0 \in C(f, \infty)$,
$w_0$ lies in a horizontal line
separating a region where $Val(f, w) = 2$
from a region where $Val(f, w) = 1$.
Thus, 
for all $w$ in a sufficiently small neighborhood 
of $w_0$,
$Val(f, w) \ge Val(f, w_0) = 1$,
in accord with the result above.

\begin{remark}
\label{rem:emptyint}
The results above are vacuous if 
the partitioning set $f(S) \cup C(f)$ fills the
plane.
When does the partitioning set have empty interior?
If  $int(f(S)) = \varnothing$,
$f(S) \cup C(f)$ will have empty interior
iff $int (C(f)) = \varnothing$.
Why?
Let $U$ be an open subset of $f(S) \cup C(f)$.
Since $C(f)$ is closed,
$U \backslash C(f)$ is an open subset of
$f(S)$.
Since $f(S)$ has empty interior,
either $U = \varnothing$ or $U$ is a non-empty
subset of $C(f)$.
\end{remark}

\begin{theorem}
\label{thm:c1dense}
Let $R$ be open in $\mathbb{R}^2$.
Let $f: R \rightarrow \mathbb{R}^2$ be a 
$C^1$ mapping such that $S$ is nowhere dense.
Suppose that $C(f)$ has non-empty interior.
Then points with infinite valence are dense in the
interior of $C(f)$.
\end{theorem}
\begin{proof}
Suppose not.
Then we can choose $w_0 \in $ int $C(f)$
and some $\epsilon_0 > 0$ such that
$f$ has finite valence at each point in 
$B(w_0, \epsilon_0)$.

Since $w_0$ is a cluster point of $f$, 
there exist $\{z_n\} \subset R$ and 
$z_0 \in \partial R \cup \{\infty\}$
such that $z_n \rightarrow z_0$ and 
$f(z_n) \rightarrow w_0$.
Since $w_0$ has finite valence, 
we may choose $z_N$ such that $w_0 \neq f(z_N) \in 
B(w_0,\frac{\epsilon_0}{2}) \subset B(w_0, \epsilon_0)$.
If $z_N \notin S$, let $\zeta_1 = z_N$.
Otherwise, 
by the continuity of $f$ in the open set $R$, 
$\exists\  
\delta > 0$
such that $B(z_N, \delta) \subset R$ and
$f(B(z_N,\delta)) \subset 
B(w_0,\frac{\epsilon_0}{2})$.
Since $f$ has finite valence at each point 
in $B(w_0, \epsilon_0)$ and since $S$
is nowhere dense in $R$, 
choose $\zeta_1 \in B(z_N,\delta) 
\backslash S$ such that $f(\zeta_1) \neq w_0$.
Let $w_1 = f(\zeta_1)$.
Since $f$ is continuous at $\zeta_1$,
we may choose $0 < \epsilon_1 < 
\frac{\epsilon_0}{2}$
and $0 < \delta_1 < 
dist(\zeta_1, \partial R)\ /\ 2$
such that 
$f(B(\zeta_1, \delta_1)) \subset 
\overline{B(w_1, \epsilon_1)}
\subset B(w_0, \epsilon_0)$.
Since $\zeta_1 \in R \backslash S$,
by the inverse function theorem, 
we can find non-empty open subsets 
$U_1 \subseteq B(\zeta_1, \delta_1)$
and $V_1 \subseteq B(w_1, \epsilon_1)$
such that $\zeta_1 \in U_1$, $w_1 \in V_1$
and $f: U_1 \rightarrow V_1$ is 1-1, onto.
Moreover, $f(U_1) = V_1 \subset
B(w_0, \epsilon_0) \subset C(f)$.

Repeat the preceding argument with $w_1$ in place
of $w_0$ to find 
$z_N \in R \backslash B(\zeta_1, \delta_1)$
such that $f(z_N) \notin \{w_0, w_1\}$ and 
$f(z_N) \in V_1$.
This gives us $\zeta_2 \in B(z_N,\delta) \backslash S$ 
such that 
$w_2 = f(\zeta_2) \notin \{w_0, w_1\}$.
Arguing as above,
we may choose $0 < \epsilon_2 < \epsilon_1$ such that 
$\overline{B(w_2, \epsilon_2)} \subset V_1$.
Similarly,
we may choose
$0 < \delta_2 < dist(\zeta_2, \partial R)\ /\ 2$
such that 
$B(\zeta_2, \delta_2) \cap B(\zeta_1, \delta_1) =
\varnothing$ and such that
$f(B(\zeta_2, \delta_2)) \subseteq B(w_2, \epsilon_2)$.
As above, we can find non-empty open subsets 
$U_2 \subseteq B(\zeta_2, \delta_2) \subset R$
and $V_2 \subseteq B(w_2, \epsilon_2)$
such that $\zeta_2 \in U_2$, $w_2 \in V_2$, and
$f: U_2 \rightarrow V_2$ is 1-1, onto.
Moreover, $f(U_2) = V_2 \subseteq B(w_2, \epsilon_2)
\subset V_1 \subset
B(w_1, \epsilon_1) \subset B(w_0, \epsilon_0) 
\subset C(f)$.
By construction,
$U_1 \cap U_2 = \varnothing$.

Continue in this manner to get a sequence of nested
non-empty open sets $V_n \subseteq B(w_n, \epsilon_n)$
with $\overline{B(w_n, \epsilon_n)} \subset 
B(w_{n-1}, \epsilon_{n-1})$
such that $w_n \rightarrow
w^* = \cap \overline{B(w_n, \epsilon_n)}
\subset B(w_0, \epsilon_0)$.
Thus, 
$Val(f, w^*)$ is finite.
But $w^*$ has a preimage in each $U_n$ 
where the $U_n$ are by construction
pairwise disjoint.  
This contradicts $w^*$ having a finite
number of distinct preimages in $R$.
\end{proof}

We are not claiming that points with
infinite valence are only in the interior of
$C(f)$.
In Example \ref{ex:flatpoly},
the origin has infinite valence and is in $C(f)$.
However, $C(f)$ has empty interior.

\subsection {Partitioning the preimage by 
$f^{-1}(f(S) \cup C(f))$}

Suppose that $R$ is an open set in $\mathbb{R}^2$
and that $f: R \rightarrow \mathbb{R}^2$ is $C^1$.
Recall that $f(S) \cup C(f)$ is closed.
We have seen that $f(S) \cup C(f)$ partitions the plane
into components of constant valence.
What does this tell us about the behavior of $f$
in a component of 
$R \backslash f^{-1}(f(S) \cup C(f))$?

\begin{theorem}
\label{thm:preimonto}
Let $R$ be an open set in $\mathbb{R}^2$.
Suppose that 
$f: R \rightarrow \mathbb{R}^2$ is a 
$C^1$ mapping.
Let $R_0$ be a connected component of 
$R\ \backslash f^{-1}(f(S) \cup C(f))$ 
and choose $z_0 \in R_0$.
Let $w_0 = f(z_0)$ and choose the connected component 
$\Omega_0 \subseteq\ 
\mathbb{R}^2\ \backslash (f(S) \cup C(f))$
such that $w_0 \in \Omega_0$.
Then $f(R_0) = \Omega_0$.
\end{theorem}
\begin{proof}
Note that the result holds vacuously if
$f^{-1}(f(S) \cup C(f)) = R$.
Suppose that $R \backslash f^{-1}(f(S) \cup C(f))
\neq \varnothing$.

We first show that $f(R_0) \subseteq \Omega_0$.
Since $R_0$ is connected and $f$ is continuous,
$f(R_0)$ is connected.
Since $w_0 \in f(R_0)$, 
by our choice of $\Omega_0$, 
$f(R_0) \subseteq \Omega_0$.

It remains to show that 
$\Omega_0 \backslash f(R_0)$ is empty.
Suppose not.
Since $R_0$ is open and $R_0 \cap S$ is empty,
$f$ is an open map on $R_0$;
hence $f(R_0)$ is open.
Since $f(R_0)$ is open and a proper subset of 
the component $\Omega_0$,
$\exists\ \tilde{w} \in \Omega_0 \cap 
(\ \overline{f(R_0)}\ \backslash\ f(R_0))$.
By Lemma \ref{lem:infval},
$Val(f, w_0)$ is finite.
Suppose that $Val(f, w_0) = N_0$.
By Theorem \ref{thm:c1part},
the valence of 
$f$ is constant in $\Omega_0$;
hence $Val(f, \tilde{w}) = N_0$.
Thus $\exists\ z_1, ..., z_{N_0}$ distinct in
$R\ \backslash f^{-1}(f(S) 
\cup C(f))$ such that $f(z_j) = \tilde{w}$ for
$j=1,..., N_0$.
Since $\tilde{w} \notin f(R_0)$,
$\{z_1, ..., z_{N_0}\} \cap R_0 = \varnothing$.

For each $z_j$, we may find an open neighborhood of 
$z_j$ in $R$, 
say $B_j$,
such that $B_j \cap S = \varnothing$,
$B_j \cap R_0 = \varnothing$,
and such that $f$ is 1-1 and an open map on $B_j$.
We may choose the $B_j$ to be pairwise disjoint.
Then $V =\ \cap_{j = 1}^{N_0} \ f(B_j)$ 
is open and non-empty.
Also, $\tilde{B_j} = f^{-1}(V) \cap B_j$
is an open neighborhood of $z_j$ where $f$ is
1-1 and such that $\tilde{B_j} \cap R_0$ is empty.
Since $\tilde{w} \in V$ and 
$\tilde{w} \in \overline{f(R_0)}\ $,
$\exists\ w \in f(R_0)$ such that
$w \in V$.
Thus, each of the pairwise disjoint $\tilde{B_j}$ 
contains one preimage of $w$.
But, $R_0$ also contains at least one preimage of
$w$.
Thus, $Val(f, w) \geq N_0 + 1$.
However, $w \in \Omega_0$, 
so $Val(f, w) = N_0$, a contradiction.
Thus, $f(R_0) = \Omega_0$.
\end{proof}

\begin{lemma}
\label{lem:n0preims}
Let $f$, $R_0$, $\Omega_0$, and $w_0$ be as in 
Theorem \ref{thm:preimonto}.
Suppose that $w_0$ has exactly $n_0$ distinct preimages
in $R_0$,
where $0 < n_0 \leq\ Val(f, w_0)$.
Then every $w \in \Omega_0$ has exactly $n_0$ distinct
preimages in $R_0$.
\end{lemma}
\begin{proof}
This proof was suggested by D. Sarason.
By Theorem \ref{thm:c1part} and Lemma \ref{lem:infval},
$Val(f, w) = Val(f,  w_0) = N_0 < \infty$
for all $w \in \Omega_0$.
By Theorem \ref{thm:preimonto},
$1 \le Val(f|_{R_0}, w) \le N_0$ for all
$w \in \Omega_0$.
Let $W_j = \{w \in \Omega_0:
Val(f|_{R_0}, w) = j\}$.
Then $\Omega_0 = \cup_{j = 1}^{N_0} W_j$.
Clearly, the $W_j$ are pairwise disjoint.
If the $W_j$ are open,
then only one of the $W_j$ is nonempty since
$\Omega_0$ is connected.
Further, 
since $Val(f|_{R_0}, w_0) = n_0$,
the result follows.

It remains to show that $W_j$ is open for $j > 0$.
This follows from Lemma \ref{lem:nvalopen}
if $W_j = V_j(f|_{R_0})$.
By construction,
the critical set of $f|_{R_0}$ is empty;
hence $V_j(f|_{R_0}) = \{w \notin C(f|_{R_0}):
Val(f|_{R_0}, w) = j\}$.
We need to show that $C(f|_{R_0})$ is disjoint from
$\Omega_0$.
Let $\{z_n\} \subset R_0$ converge to a point
$z_0 \in \partial R_0 \cup \{\infty\}$.
If $z_0 \in \partial R \cup \{\infty\}$,
then if $\{f|_{R_0}(z_n)\}$ has a finite cluster point,
this point is in $C(f)$,
hence not in $\Omega_0$.
On the other hand,
if $z_0 \in R$,
then $f$ is continuous at 
$z_0 \in f^{-1}(f(S) \cup C(f))$
and $\{f|_{R_0}(z_n)\}$ converges to a point
in $f(S) \cup C(f)$,
which, again, is not in $\Omega_0$.
Thus $V_j(f|_{R_0}) = W_j$ and $W_j$ is open.
\end{proof}

Thus,
$f^{-1}(f(S) \cup C(f))$ partitions $R$ into
components,
each of which is mapped onto a component
of $\mathbb{R}^2 \backslash (f(S) \cup C(f))$
by $f$.
We will now show that $f$ is a covering map
on each component of this partition of the preimage.
First,
we recall a standard result
(see~\cite{Mun:top}, page 341):

\begin{theorem}[Munkres]
\label{thm:munch8th4.5}
Let $p: (E, e_0) \rightarrow (B, b_0)$ be a 
covering map.
If $E$ is path connected,
then there is a surjection
$\phi: \pi_1(B, b_0) \rightarrow p^{-1}(b_0)$.
If $E$ is simply connected,
$\phi$ is a bijection.
\end{theorem}

\begin{theorem}
\label{thm:preimcover}
Let $R$ be an open set in $\mathbb{R}^2$.
Suppose that $f: R \rightarrow \mathbb{R}^2$
is $C^1$.
Let $R_0$ be a connected component of 
$R \backslash f^{-1}(f(S) \cup C(f))$.
Choose $z_0 \in R_0$ and let
$w_0 = f(z_0)$.
Suppose that $w_0$ has exactly $n_0$ distinct
preimages in $R_0$.
Then $R_0$ is a $n_0$-fold covering of $\Omega_0$.
Moveover, if $\Omega_0$ is simply connected,
then $f$ is univalent in $R_0$.
\end{theorem}
\begin{proof}
We first show that $f: R_0 \rightarrow \Omega_0$
is a covering map.
By Theorem \ref{thm:preimonto}, $f$ is onto;
$f$ is continuous by assumption.
Choose $w_0 \in \Omega_0$.
By Lemma \ref{lem:n0preims}, 
$w_0$ has $n_0$ distinct preimages in 
$R_0$.
Constructing $V$ and the $\tilde{B_j}$ as in 
the proof of Theorem \ref{thm:preimonto}
(except that we choose $B_j \subset R_0$),
each of the pairwise disjoint open sets $\tilde{B_j}$ 
is homeomorphic to $V$
since $R_0 \cap S = \varnothing$.
Also, by Lemma \ref{lem:n0preims},
each $w \in V$ has exactly $n_0$ distinct preimages
in $R_0$.
Thus, $R_0 \cap f^{-1}(V) = 
\cup_{j = 1}^{n_0} \tilde{B_j}$.
Thus, $V$ is evenly covered by $f|_{R_0}$.
Since $w_0$ is arbitrary,
$f|_{R_0}$ is a covering map
and $R_0$ is a covering space of $\Omega_0$.
Since $f^{-1}(w)$ has $n_0$ distinct elements in
$R_0$ for each $w \in \Omega_0$,
$R_0$ is a $n_0$-fold covering of $\Omega_0$.

We note that $R \backslash f^{-1}(f(S) \cup C(f))$
is open in $\mathbb{R}^2$.
Hence the component $R_0$ is open in $\mathbb{R}^2$.
Since $R_0$ is an open, connected subset of 
$\mathbb{R}^2$,
given any two points in $R_0$,
we may find a polygonal path contained in $R_0$
that joins the two points.
Hence, $R_0$ is path connected.
Choose $w_0 \in \Omega_0$. 
By Theorem \ref{thm:munch8th4.5}, 
there is a surjection
$\phi: \pi_1 (\Omega_0, w_0) \rightarrow f^{-1}(w_0)$,
where we are restricting $f$ to $R_0$.
If we also assume that $\Omega_0$ is simply
connected,
then $\pi_1 (\Omega_0, w_0)$ is the trivial group.
Since $\phi$ is a surjection,
$w_0$ must have exactly one preimage in $R_0$.
Hence $n_0 = 1$ and $f$ is univalent on $R_0$.
\end{proof}

\begin{example}
\label{ex:c1trans}
$f(x,y) = (x \cos y,\ y)$

Here,
$R = \mathbb{R}^2$ and $C(f) = C(f, \infty)$.
If we rewrite $f$ as $f = u + iv$ where $z = x + iy$, 
it is clear that $f(z)$ is $C^1$ but not harmonic.
A calculation shows that
\[
\begin{array}{l}
\ S = \{(x, \frac{(2k + 1) \pi}{2}): 
x \in \mathbb{R}$ and $
k \in \mathbb{Z}\}\\
\ f(S) = \{(0, \frac{(2k + 1) \pi}{2}):
k \in \mathbb{Z}\}\\
\end{array}
\]
\noindent
Also,
$C(f, \infty)$ consists of the horizontal lines
$y = \frac{(2k + 1) \pi}{2}$ 
where $k$ is an integer.
Each 
$w \in \mathbb{R}^2 
\backslash (f(S) \cup C(f, \infty))$
has exactly one preimage.
Each $w \in f(S)$ has an infinite number of preimages.
Each $w \in C(f, \infty) \backslash f(S)$ has no
preimages.
Note that if we fix $a \in \mathbb{R}$  
and let $y_n = \cos^{-1}(a / x_n)$,
then $f(x_n, y_n) \rightarrow (a, 
\frac{(2k + 1) \pi}{2}) \in C(f, \infty)$
as $x_n \rightarrow \infty$,
provided that we
choose the branch of $\cos^{-1}$ such that
$\cos^{-1}(a / x_n) \rightarrow 
\frac{(2k + 1) \pi}{2}$.

The partitioning set of our domain is
$f^{-1}(f(S) \cup C(f, \infty))$,
which is $S$ (a collection of horizontal lines.)
It is clear that if $(x, y) \notin S$,
then $f(x, y) \neq (0, \frac{2 k + 1}{2} \pi)$.
If we choose $w = (a, b) \in \mathbb{R}^2$
such that $w \notin f(S) \cup C(f, \infty)$,
then $b \neq \frac{2 k + 1}{2} \pi$.
It is clear that $w$ has exactly one preimage;
namely,
$(a / \cos b, b)$ and that this preimage point
does not lie in $S$.
Hence,
$f^{-1}(f(S) \cup C(f, \infty))$ partitions
$\mathbb{R}^2$ into horizontal strips where
$f$ is univalent.
\end{example}

\subsection {Adjacent components of the preimage}
A non-empty, connected set is said to be degenerate
if it consists of a single point.
If two distinct components of our partition of $R$
share a non-degenerate common boundary arc in $R$,
will $f$ map both to the same component of 
$\mathbb{R}^2 \backslash (f(S) \cup C(f))$?

\begin{lemma}
\label{lem:preimadj}
Let $R \subseteq \mathbb{R}^2$ be open.
Let $f: R \rightarrow \mathbb{R}^2$ be $C^1$. 
Let $R_1$ and $R_2$ be distinct components of
$R\ \backslash\ 
f^{-1} (f(S) \cup C(f))$
such that
$\overline{R_1} \cap \overline{R_2} \neq \varnothing$.
Suppose that $f(R_1) = f(R_2) = \Omega$.
If $\partial \Omega \cap int\ \overline{\Omega}
= \varnothing$,
then there exists no non-empty set 
$\gamma \subseteq\ 
(\overline{R_1}\cap\overline{R_2} \cap R)\ 
\backslash\ S$
such that
$R_1 \cup R_2 \cup \gamma$ is open.
\end{lemma}
\begin{proof}
By contradiction.
Suppose that $R_1$ and $R_2$ are disjoint components
of $R\ \backslash\ 
f^{-1} (f(S) \cup C(f))$
such that  
$f(R_1) = f(R_2)$.
Suppose also that there exists 
$\gamma \subseteq
(\overline{R_1}\cap\overline{R_2} \cap R)\  
\backslash\ S$
such that $\gamma \neq \varnothing$ and
such that $R_1 \cup R_2 \cup \gamma$ is open.
By Theorem \ref{thm:preimonto}, there exist
$\Omega_1, \Omega_2$, components of
$\mathbb{R}^2\ \backslash\ (f(S) \cup C(f))$,
such that $f(R_1)=\Omega_1$ and $f(R_2)=\Omega_2$.
By assumption, $\Omega_1=\Omega_2 = \Omega$.

Since $R_1$ and $R_2$ are disjoint components of
$R\ \backslash\ 
f^{-1} (f(S) \cup C(f))$
and $\gamma \subseteq \overline{R_1} \cap 
\overline{R_2} \cap R$,
$\gamma \subseteq f^{-1} (f(S) \cup C(f))$.
Thus $f(\gamma) \subseteq f(S) \cup C(f)$.
By continuity, $f(\gamma) \subseteq \overline{\Omega}$,
so $f(\gamma) \subseteq \partial \Omega$.
Now let $R_0 = R_1 \cup R_2 \cup \gamma$.
By assumption,
$R_0$ is open.
Since $R_0 \cap S = \varnothing$,
$f$ is an open map on $R_0$ by the inverse function
theorem and $f(R_0)$ is open.
But $f(R_0) = f(R_1) \cup f(R_2) \cup f(\gamma) 
= \Omega \cup f(\gamma)$
and we have seen that 
$f(\gamma) \subseteq \partial \Omega$.
Since $f(R_0)$ is open,
$f(R_0) \subseteq int\ \overline{\Omega}$.
Choose $z_0 \in \gamma$
(recall that $\gamma \neq \varnothing$.)
Then $f(z_0) \in \partial \Omega \cap f(R_0)
\subseteq \partial \Omega \cap int\ 
\overline{\Omega}$.
Thus $\partial \Omega \cap int\ \overline{\Omega} 
\neq \varnothing$,
a contradiction.
\end{proof}

\begin{lemma}
\label{lem:harmadj}
Let $R \subseteq \mathbb{R}^2$ be open.
Let $f: R \rightarrow \mathbb{R}^2$ be a 
light $C^1$ function.
Let $R_1$ and $R_2$ be distinct components of
$R\ \backslash\ 
f^{-1} (f(S) \cup C(f))$
such that
$\overline{R_1} \cap \overline{R_2} \neq \varnothing$.
Suppose that $f(R_1) = f(R_2) = \Omega$.
If $\partial \Omega \cap int\ \overline{\Omega}$
consists of a finite number of points,
then there exists no non-empty,
non-degenerate connected set 
$\gamma \subseteq\ 
(\overline{R_1}\cap\overline{R_2} \cap R)\ 
\backslash\ S$
such that
$R_1 \cup R_2 \cup \gamma$ is open.
\end{lemma}
\begin{proof}
By contradiction.
Suppose that $R_1$ and $R_2$ are disjoint components
of $R\ \backslash\ 
f^{-1} (f(S) \cup C(f))$
such that $f(R_1) = f(R_2)$.
Suppose also that there exists 
$\gamma \subseteq
(\overline{R_1}\cap\overline{R_2} \cap R)\ 
\backslash\ S$
such that $\gamma \neq \varnothing$ and
such that $R_0 = R_1 \cup R_2 \cup \gamma$ is open.

Let $P = \partial \Omega \cap int\ \overline{\Omega}$.
By assumption, $P = \{w_1, w_2, ... w_n\}$.
If $n = 0$, the result follows from 
Lemma \ref{lem:preimadj}.
Assume that $n > 0$.
By assumption, $\gamma$ is a non-empty, 
non-degenerate connected set.
Since $f$ is continuous, $f(\gamma)$ is connected.
So, if $f(\gamma) \subseteq P$,
then $f(\gamma) = \{w_j\}$ for some fixed value of $j$.
This contradicts $f$ being a light mapping.
Thus $f(\gamma) \backslash P \neq \varnothing$.
Choose $w \in f(\gamma) \backslash P$.
By the arguments used in the proof of 
Lemma \ref{lem:preimadj},
$f(\gamma) \subseteq \partial \Omega$.
As in Lemma \ref{lem:preimadj}, $f(R_0)$ is open and 
$w \in int\ \overline{\Omega}$.
Thus $w \in \partial \Omega \cap 
int\ \overline{\Omega} = P$,
a contradiction.
\end{proof}

\begin{corollary}
\label{cor:planadj}
Let $f$ be a light $C^1$ function in $\mathbb{R}^2$.
Let $\Omega$ be a component of
$\mathbb{R}^2 \backslash (f(S) \cup C(f, \infty))$.
Suppose that $R_1$ and $R_2$ are two
distinct components of 
$\mathbb{R}^2 \backslash 
f^{-1}(f(S) \cup C(f, \infty))$
such that $f(R_1) = f(R_2) = \Omega$.
\begin{enumerate}
\item
If $\partial \Omega \cap int\ \overline{\Omega}$
consists of a finite number of points,
then $(\overline{R_1} \cap \overline{R_2}) 
\backslash S$
contains no non-degenerate connected set $\gamma$
such that $R_1 \cup R_2 \cup \gamma$ is open.
\item
Suppose that
$\partial \Omega \cap int\ \overline{\Omega}$
consists of a finite number of points.
Also suppose that
$(int\ \overline{R_j}) \cap \partial \overline{R_j}$
consists of a finite number of points and that
$\partial (int\ \overline{R_j})$ is a Jordan curve
for $j = 1, 2$.
If $\overline{R_1} \cap \overline{R_2}$
contains a non-degenerate Jordan arc,
then this arc is in $S$.
\end{enumerate}
\end{corollary}
\begin{proof}
Using the notation of the preceding two lemmas,
$R = \mathbb{R}^2$.
So $\overline{R_1} \cap \overline{R_2} \cap R = 
\overline{R_1} \cap \overline{R_2}$.
The first claim follows from Lemma \ref{lem:harmadj}.

The second claim follows by noting that each
region $R_j$ is a Jordan region that has a finite number
of puncture points.
We are assuming that the boundaries intersect 
in a non-degenerate Jordan arc $\rho$.
Since we can find subregions of the two regions,
each having $\rho$ as part of its boundary,
that are Jordan regions (no punctures),
$R_1 \cup R_2 \cup int\ \rho$ is open.
By the first claim,
we must have that $(int\ \rho) \backslash S$
is totally disconnected.
Since $S$ is closed, $\rho \subseteq S$.
\end{proof}

\subsection{Behavior in shared boundaries and
``puncture points''}

Let $f$ be a $C^1$ function on some open
set $R \subseteq \mathbb{R}^2$.
Suppose that we have a collection of components
of $R \backslash f^{-1}(f(S) \cup C(f))$
where $f$ is univalent on each component.
Can we combine these components to get a larger region
where $f$ is univalent?
The results in the preceding section 
(such as Corollary \ref{cor:planadj})
handle components that share a common boundary arc
in $R$,
but say nothing about the behavior on the boundary arc. 

If $R_0$ is one such component,
we want to know more about the behavior of
$f$ at points in $\partial R_0 \cap R$.
These could be points in boundary arcs shared
with other components or could be ``puncture points.''
A ``puncture point''
is an isolated point in $\partial R_0$.
Does $f$ map $\partial R_0$ onto $\partial(f(R_0))$?
What can we say about $Val(f|_{\partial R_0}, w)$
for $w \in \partial(f(R_0))$?

\begin{lemma}
\label{lem:harmbdyim0}
Let $R \subseteq \mathbb{R}^2$ be open.
Let $f$ be $C^1$ in $R$.
Let $R_0 \subset \mathbb{R}^2$ be a component of
$R\ \backslash\ 
f^{-1} (f(S) \cup C(f))$.
Then $f(R \cap \partial R_0) \subseteq 
f(R) \cap \partial(f(R_0))$.
\end{lemma}
\begin{proof}
By Theorem \ref{thm:preimonto}, 
there exists a component of
$\mathbb{R}^2\ \backslash\ (f(S) \cup C(f))$,
say $\Omega_0$,
such that $f(R_0) = \Omega_0$.
Thus $\partial (f(R_0)) = \partial \Omega_0$.
Requiring that $R_0$ is a proper subset of 
$\mathbb{R}^2$
guarantees that $\partial R_0 \neq \varnothing$
in the finite plane.
If $\partial R_0 \subseteq \partial R$,
the result holds vacuously.
Suppose that $R \cap \partial R_0 \neq \varnothing$.

Since $R \cap \partial R_0 \subseteq
f^{-1} (f(S) \cup C(f))$,
$f(R \cap \partial R_0) \subseteq f(S) \cup C(f)$.
Fix $z_0 \in R \cap \partial R_0$.
Then there exists a sequence $\{z_n\} \subset R_0$
such that $z_n \rightarrow z_0$.
Since $z_0 \in R$,
by the continuity of $f$ at $z_0$,
$f(z_n) \rightarrow f(z_0)$.
But $f(z_n) \in \Omega_0$ for all $n$ and
$f(z_0) \in f(S) \cup C(f)$.
Hence $f(z_0) \in 
f(R) \cap \overline{\Omega_0} \cap
(f(S) \cup C(f))\ \subseteq\ 
f(R) \cap \partial \Omega_0$,
so $f(R \cap \partial R_0) \subseteq\ 
f(R) \cap \partial \Omega_0$.
\end{proof}

\begin{lemma}
\label{lem:harmbdyim}
Let $f$ be $C^1$ in $\mathbb{R}^2$.
Let $R_0$ be a bounded component of
$\mathbb{R}^2 \backslash 
f^{-1} (f(S)\\ \cup C(f, \infty))$.
Then $f(\partial R_0) = \partial(f(R_0))$.
\end{lemma}
\begin{proof}
By Theorem \ref{thm:preimonto}, 
there exists a component of
$\mathbb{R}^2\ \backslash\ (f(S) \cup C(f, \infty))$,
say $\Omega_0$,
such that $f(R_0) = \Omega_0$.
Thus $\partial (f(R_0)) = \partial \Omega_0$.
By Lemma \ref{lem:harmbdyim0},
we have $f(\partial R_0) \subseteq
f(\mathbb{R}^2) \cap \partial \Omega_0$.
So $f(\partial R_0) \subseteq \partial \Omega_0$.

We now show that 
$\partial \Omega_0 \subseteq f(\partial R_0)$.
Fix $w_0 \in \partial \Omega_0$.
Then there exists a sequence $\{w_n\} \subset \Omega_0$
such that $w_n \rightarrow w_0$.
Since $f(R_0) = \Omega_0$,
there exist $z_n \in R_0$ such that
$f(z_n) = w_n$ for each $n$.
By assumption $R_0$ is bounded,
so $\overline{R_0}$ is compact and $\{z_n\}$ has
a convergent subsequence $\{z_{n_k}\}$.
Since $w_n \rightarrow w_0$,
$w_{n_k} \rightarrow w_0$.
Let $z_{n_k} \rightarrow\ z_0 \in \overline{R_0}$.
Since $f$ is continuous in $\mathbb{R}^2$,
$f(z_0) = w_0$.
Since $w_0 \in f(S) \cup C(f, \infty)$,
$z_0 \in f^{-1} (f(S) \cup C(f, \infty))$.
Thus $z_0 \in \overline{R_0} \cap 
(f^{-1} (f(S) \cup C(f, \infty)))\ =\ \partial R_0$.
Since $w_0 \in \partial \Omega_0$ is arbitrary,
we have
$\partial \Omega_0 \subseteq f(\partial R_0)$.
The result follows.
\end{proof}

\noindent
\begin{notation}
Let $\Omega \subseteq \mathbb{R}^2$
and $f$ be a function defined on 
$R \subseteq \mathbb{R}^2$.
Recall that $Val(f|_R, \Omega)$
denotes the valence of $f|_R$ in $\Omega$.
Note that if the valence of $f|_R$ is constant and 
finite in $\Omega$,
then $Val(f|_R, \Omega) = Val(f|_R, w)$ for any
choice of $w \in \Omega$.
\end{notation}

\begin{theorem}
\label{thm:harmbdypreim}
Let $f$ be $C^1$ in $\mathbb{R}^2 $.
Let $R$ be a component of $\ \mathbb{R}^2\ \backslash
f^{-1}(f(S) \cup C(f, \infty))$ and let
$\Omega$ be a component of $\mathbb{R}^2\ \backslash
(f(S) \cup C(f,\infty))$
such that $f(z) \in \Omega$ for some $z \in R$.
Choose $w \in \partial \Omega$ such that
\emph{(i)} $w \notin C(f\mid_{R},\ \infty)$,
\emph{(ii)} $f^{-1}(w)\ \cap\ \partial R\ \cap\ S 
= \varnothing$,
\emph{(iii)} $f^{-1}(w)\ \cap\ \partial R 
\neq \varnothing$,
and such that
\emph{(iv)} for each $\zeta \in  f^{-1}(w) \cap 
\partial R$,
there exists a neighborhood of $\zeta$, 
say $U_\zeta$,
such that
$f(U_\zeta\ \backslash\ R)\ \cap\ \Omega = \varnothing$.
Then $w$ has exactly $Val(f\mid_{R}, \Omega)$ 
distinct preimages on $\partial R$.
\end{theorem}
\begin{proof}
By Theorem \ref{thm:preimonto}, $f(R) = \Omega$. 
Moreover, by Lemma \ref{lem:n0preims},
if $\eta \in \Omega$ and if $N = Val(f\mid_{R}, \eta)$, 
then $Val(f\mid_{R}, \Omega) = N$.
Choose $w \in \partial \Omega$ satisfying the
given conditions.
By condition (iii),
$w$ has at least one preimage on $\partial R$.
Note that if $R$ is bounded,
this condition follows from 
Lemma \ref{lem:harmbdyim}.
Suppose that $w$ has $M$ distinct preimages on 
$\partial R$.

We first show that $M \ge N$:
Since $w \in \partial \Omega$,
we may choose $\{w_j\} \subset \Omega$,
where the $w_j$ are distinct,
such that $w_j \rightarrow w$.
For each $j$,
there exist $z_{j1}, z_{j2}, \ldots , z_{jN}$ distinct
in $R$ such that $f(z_{jk}) = w_j$.
Let $Q = \bigcup_{k=1}^{N} 
\{ z_{jk} \}_{j=1}^{\infty}$.
So $Q \subset R$.
Using arguments analogous to those
used to prove Lemma \ref{lem:nvalopen},
we see that $Q$ has at least $N$ distinct
cluster points.
We also see that the cluster points of $Q$
are in $\partial R$.
Hence
$M \ge N$.

We now show that $M \le N$:
Suppose not.
Then we can find distinct values $z_1, \ldots, z_{N+1}$
such that $z_j \in \partial R$ and $f(z_j) = w$.
By (ii),
$z_j \notin S$.
So, for each $j$,
we may find some neighborhood of $z_j$, 
say $V_j$,
where $f$ is a homeomorphism.
We may choose the $V_j$ so that they are pairwise
disjoint.
By (iv),
for each $z_j$,
there exists a neighborhood $U_j$ of $z_j$ such
that the preimages in $U_j$ of all points in 
$f(U_j)\ \cap\ \Omega$
lie in $R$.
Let $B_j = U_j\ \cap\ V_j$.
Then the $B_j$ are pairwise disjoint and each 
$B_j$ is a non-empty, open neighborhood
of $z_j$ where $f$ is a homeomorphism.
Since $z_j \in \partial R$,
$B_j\ \cap R\ \neq\ \varnothing$;
hence
$f(B_j)\ \cap\ \Omega\ \neq\ \varnothing$.
By construction, 
all preimages in $B_j$ of points in
$f(B_j)\ \cap\ \Omega$ lie in $R$.
Since $z_j \in B_j$ and $f(z_j) = w$,
$\bigcap_{j = 0}^{N+1} f(B_j) \neq \varnothing$.
Choose $\eta \in \Omega\ \cap\ 
(\bigcap_{j = 1}^{N + 1} f(B_j))$.
Then $\eta$ has a distinct preimage in each $B_j$,
giving a total of $N+1$ distinct preimages in 
$\bigcup_{j=1}^{N+1} B_j$.
Further, 
each of these preimages lies in $R$.
But,
$\eta$ has exactly $N$ distinct preimages in $R$,
a contradiction.
Thus, $M = N$ and $w$ has exactly 
$Val(f\mid_{R}, \Omega)$
distinct preimages on $\partial R$.
\end{proof}

\noindent
\begin{remark}
We note that
\begin{enumerate}
\item
We will apply this result to prove
Lemma \ref{lem:puncturepoints} below.
\item
The condition $w \notin C(f\mid_{R},\ \infty)$ 
automatically holds if either $C(f, \infty)$ is empty
or if $R$ is bounded.
\end{enumerate}
\end{remark}

\begin{theorem}
\label{thm:jointwoc1}
Let $f$ be a light, $C^1$ function in $\mathbb{R}^2$.
Let $R_1$ and $R_2$ be disjoint bounded
components of 
$\ \mathbb{R}^2 \backslash 
f^{-1}(f(S) \cup C(f, \infty))$.
Suppose that $R_1$ and $R_2$ each have at most a
finite number of puncture points.
Suppose that for each component $\Omega \subseteq
\mathbb{R}^2 \backslash (f(S) \cup C(f, \infty))$,
$\partial \Omega \cap int\ \overline{\Omega}$ 
consists of a finite number of points.
Suppose that $f$ is univalent on $R_1$ and on $R_2$.
Let $\gamma$ be a non-degenerate, simple arc
in $\partial R_1 \cap \partial R_2$.
Suppose that $f(int\ \gamma)$ is a simple arc.
If $S\ \cap\ int\ \gamma = \varnothing$,
then
\begin{enumerate}
\item
$f(R_1)$ and $f(R_2)$ are distinct components
of $\ \mathbb{R}^2 \backslash (f(S) \cup C(f,\infty))$.
\item
$f$ is a univalent on 
$R_1\ \cup\ R_2\ \cup\ int\ \gamma$. 
\end{enumerate}
\end{theorem}
\begin{proof}
By construction, the combined region, 
$R_1\ \cup\ R_2\ \cup\ int\ \gamma$,
is open
(see the proof of Corollary \ref{cor:planadj}.)
The first claim follows from 
Lemma \ref{lem:harmadj}.
To prove the second claim,
we need to show that $f$ is univalent on $int\ \gamma$.
Then $f$ is univalent in the combined region
(since no point in $R_1 \cup R_2$ is mapped to
$f(int\ \gamma)$.)
Note that $f(int\ \gamma)$ is non-degenerate
since $f$ is light.

We now check that $f$ is univalent on $int\ \gamma$.
Suppose not.
Then we have two distinct points $\zeta_1,
\zeta_2 \in int\ \gamma$ such that
$f(\zeta_1) = f(\zeta_2) = w_0 
\in f(S) \cup C(f, \infty)$.
Since $S \cap\ int\ \gamma = \varnothing$,
we can find a neighborhood $B_1$ of $\zeta_1$
and a neighborhood $B_2$ of $\zeta_2$
where $f$ is an open map.
We may assume that these neighborhoods are disjoint
and are contained in the combined region.
However,
$f(B_1) \cap f(B_2)$ is an open neighborhood
of $w_0$.
So it contains a point $w \in \Omega_1$.
By construction,
$w$ has at least two distinct preimages in $R_1$
(one in $B_1 \cap R_1$ and another in
$B_2 \cap R_1$),
which contradicts the univalence of $f$ in $R_1$.
\end{proof}

\subsection{Extension of partitioning results 
to $C^1$ mappings in $\mathbb{R}^n$}

The notation defined at the beginning
of Section \ref{sec:notation} for the critical
set and cluster set has obvious generalizations
to the case of $C^1$ functions defined in an 
open set in $\mathbb{R}^n$,
where $n$ is a positive integer.
Moreover, 
the comments at the beginning of 
Section \ref{sec:notation} concerning the properties of
these sets hold.
For the most part,
the proofs in Section \ref{sec:c1part} do not
make use of special properties of $\mathbb{R}^2$.
Thus,
the partitioning results above can be 
extended to $C^1$ functions defined in open
subsets of $\mathbb{R}^n$, 
with the exceptions of Lemma \ref{lem:harmadj},
Corollary \ref{cor:planadj},
and Theorem \ref{thm:jointwoc1}.

Lemma \ref{lem:harmadj} and the first part of
Corollary \ref{cor:planadj} hold in 
$\mathbb{R}^n$ for $n > 1$.
The hypothesis concerning the non-degenerate set
$\gamma$ cannot hold for $n = 1$.
The proofs of Theorem \ref{thm:jointwoc1}
and of the second part of Corollary \ref{cor:planadj} 
use Jordan curves and Jordan regions;
these are $n = 2$ results.

\section {Partitioning results for
harmonic mappings in an open set}
\label{sec:harm}
We will apply our results about $C^1$ functions
in an open subset of the real plane to the case
where our function is harmonic in an open subset
of the complex plane.
Unless explicitly stated otherwise,
a harmonic function is a complex-valued
harmonic function.

\begin{theorem}
\label{thm:harmpart}
Let $f$ be a harmonic mapping in an open set 
$R \subseteq \mathbb{C}$.
Then $ f(S) \cup C(f) $ partitions $\mathbb{C}$ 
into regions of constant valence.
\end{theorem}
\begin{proof}
Identify $\mathbb{R}^2$ with $\mathbb{C}$ by putting 
$x = Re\ z$ and $y = Im\ z$.  
Let $u = Re\ f$, $v = Im\ f$.
Then apply Theorem \ref{thm:c1part}.
\end{proof}

\noindent
Analogs to the other results in the preceding 
section are shown similarly.

\begin{example} 
\label{ex:transharm}
A transcendental harmonic function 
with finite valence.
\[f(z)= z + Re\ e^z\]

This example was given by D. Bshouty, W. Hengartner and 
M. Na\-ghi\-bi-Bei\-dokh\-ti
in\-~\cite{BHN:dilat}
as an example of a 2-valent,
transcendental harmonic mapping in $\mathbb{C}$. 
Calculations show that
\[C(f,\infty)=\{\zeta : Im\ \zeta = \frac{(2k + 1)\pi}
{2},\ k \in \mathbb{Z}\}\]
and
\[S = \{z : Re\ e^z = -1\}\]
Letting $z=x+iy$ and fixing $\zeta=a+ib$, 
$f(z)=\zeta$ when $y=b$ and $x + e^x \cos b=a$. 
The behavior of $\varphi(x) = x + e^x \cos b$ 
depends on the sign of $\cos b$ and determines the
number of preimages of $\zeta$.
One can then check that
$f(S) \cup C(f,\infty)$ partitions $\mathbb{C}$
into regions of constant valence.
In particular,
each $w \in f(S) \cup C(f, \infty)$ has exactly 
one preimage.
If $Im(w)\ mod\ 2\pi \in (-\pi/2, \pi/2)$, 
$w$ also has exactly one preimage.
If $Im(w)\ mod\ 2\pi \in (\pi/2, 3\pi/2)$,
then $w$ has two distinct preimages if $w$ lies
to the left of $f(S)$ and no preimages if $w$ lies
to the right of $f(S)$.
Note that if $\zeta = a + i\ \frac{(2 k + 1) \pi}{2}
\in C(f, \infty)$,
then $f(z_n) \rightarrow \zeta$ for
$z_n = x_n + i\ \cos^{-1}((a  -x_n) \exp(-x_n))$
with $x_n \rightarrow \infty$ and the proper choice
for the branch of $\cos^{-1}$.
\end{example}

\begin{example} 
\label{ex:polyunbdds}
A finite valence harmonic 
polynomial with an unbounded critical set.
\[f(z) = 2(Re(z^3) + iz)\]

One can explicity check that 
\[\{\zeta: Im\ \zeta = 0\} \cup 
\{\zeta: Re\ \zeta = \frac{(Im\ \zeta)^3}{4} + 
\frac{1}{3\ Im\ \zeta}\ \}\]
partitions $\mathbb{C}$ into regions of constant
valence.

One can also check that  
$S = \{z: xy = -\frac{1}{6},\ $where$\ x = Re\ z\ $
and $y = Im\ z\}$
and that $C(f,\infty) = \{\zeta: Im\ \zeta = 0\}$.
For example,
fix $a \in \mathbb{R}$;
then $f(z_n) \rightarrow a \in C(f, \infty)$
if $z_n = x_n\ +\ i\ \{a / 2\ -\ (1/ (3 x_n))\}$
and $x_n \rightarrow 0$.

A calculation shows that
$f(S) = \{2 (Re(z))^3 + \frac{1}{6Re(z)} + 2i Re(z): 
Re(z) \neq 0\}$,
which is equivalent to 
$\{\zeta: Re\ \zeta = \frac{(Im\ \zeta)^3}{4} + 
\frac{1}{3\ Im\ \zeta}\ \}$.
Thus $f(S) \cup C(f\infty)$ partitions $\mathbb{C}$
into regions of constant valence.
Each $w \in C(f, \infty)$ has exactly one preimage,
as does each $w \in f(S)$.
Looking at $f(S)$ in $\mathbb{R}^2$,
one sees that it consists of two parabola-like
curves pointing sideways.
Each $w$ in the region ``inside" these curves has
no preimage.
Outside of this region, each $w \notin (f(S) \cup 
C(f, \infty))$ has two distinct preimages.

Since $C(f,\infty) \neq \varnothing$, 
$f$ is an example of a finite valence harmonic
polynomial where $f$ does not go to $\infty$ as
$z \rightarrow \infty$.

It is easy to check that $f^{-1}(f(S)) = S$.
Another easy calculation shows that
$f^{-1}(C(f, \infty))$ is the imaginary axis.
Recall that $w \notin f(S) \cup C(f, \infty)$
has either zero or exactly two distinct preimages.
Suppose we choose $w \notin f(S) \cup C(f, \infty)$
such that $Val(f, w) = 2$.
Note that $Im\ w \neq 0$, 
since $w \notin C(f, \infty)$.
Since $f(x, y) = 2 (x^3 - 3 x y^2 - y + i x)$,
$Re(f^{-1}(w)) = \frac{1}{2}\ Im\  w$.
Since $w \notin f(S)$,
$Im(f^{-1}(w)) \neq -1\ / (3\ Im\ w)$.
The vertical line $Re\ z = \frac{1}{2}\ Im\ w$
intersects $S$ at 
$(\frac{1}{2}\ Im\ w, -1 / (3\ Im\ w))$.
The two preimages of $w$ are the endpoints of 
a vertical line segment bisected by $S$. 
From this,
we see that $f^{-1}(f(S) \cup C(f, \infty))$
partitions $\mathbb{C}$ into regions where
$f$ is univalent.
\end{example}

\begin{example}
\label{ex:flatpoly}
$f(z) = 2[Re(z^2) + i (Re(z) - Im(z))]$

One can readily check that
$S = \{z: Re(z) = Im(z)\}$ and $f(S) = \{0\}$. 
Further, $C(f, \infty)$ is the real axis.
If we rewrite $f$ letting $z = x + iy$,
we see that
$f(x, y) = 2 (x - y) (x + y + i)$.
A calculation shows that if $w \notin 
\mathbb{R}$, then $w$ has exactly one preimage.
If $w$ is a non-zero real, it has no preimage.
And, the origin is a point of infinite valence.
Also, if we fix $a \in \mathbb{R}$ and let 
$y_n = x_n\ -\ (a / (4 x_n))$,
then $f(z_n) \rightarrow a \in C(f, \infty)$
for $z_n = x_n + i\,y_n$ as 
$|x_n| \rightarrow \infty$.
\end{example}

\begin{figure}
\includegraphics[scale=.95]{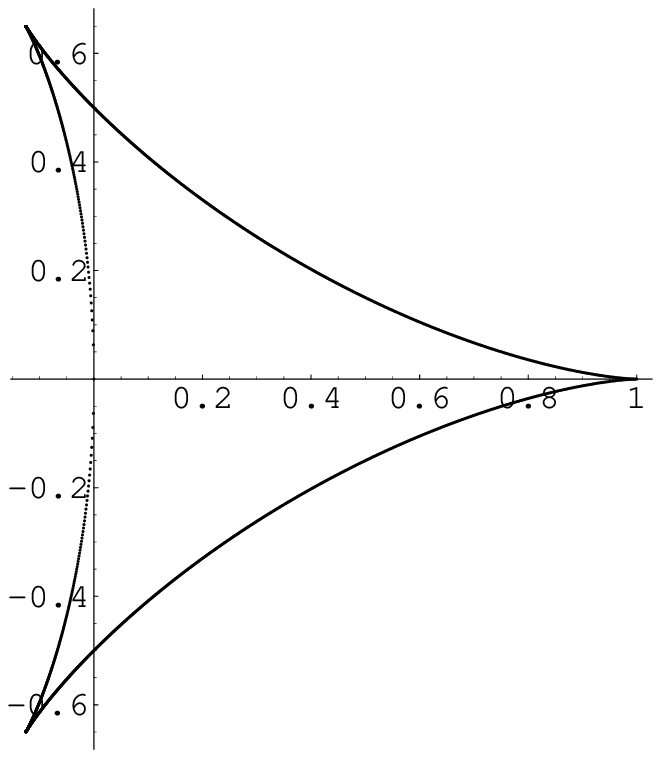}
\caption{Image of critical set for $f(z) = 
z^2 + z - \overline{z}$}
\label{fig:quadratic_imag}
\end{figure}

\begin{figure}
\includegraphics[scale=.95]{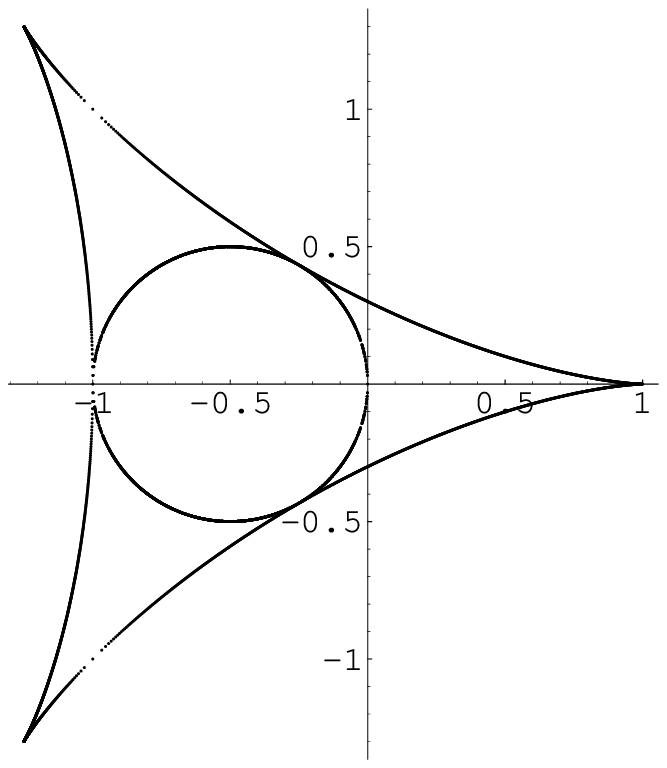}
\caption{Preimage of $f(S)$ for $f(z) = 
z^2 + z - \overline{z}$}
\label{fig:quadratic_preim}
\end{figure}

\begin{example}
\label{ex:quadraticpreim}
$f(z) = z^2 + z - \overline{z}$

An easy calculation shows that 
$S = \{z: |z + \frac{1}{2}| = \frac{1}{2}\}$.
This circle is mapped to the hypocycloid
in Figure \ref{fig:quadratic_imag}.
Clearly,
$C(f, \infty) = \varnothing$.

D. Sarason (personal communication)
has analyzed the behavior of $f$.
Among many observations about $f$,
he shows that $f$ is 1-1 on $S$,
4-valent in the bounded component of 
$\mathbb{C} \backslash f(S)$,
and 2-valent on the portion of the real axis
in the unbounded component.
Now consider our partition of the domain
(by $f^{-1}(f(S))$, since $C(f, \infty) = \varnothing$);
this is shown in Figure \ref{fig:quadratic_preim}.
Sarason's analysis shows that $f$ is a 1-1 mapping of
the interior of each of the bounded components onto
the bounded component of $\mathbb{C} \backslash f(S)$;
hence $f$ is univalent on each of the bounded
components of the partition of the domain.
He also showed that all $g(z) = p(z) - \overline{z}$,
where $deg\ p = 2$,
can be reduced to $f(z)$ by means of affine 
transformations of the domain and range.

By Theorem \ref{thm:harmpart},
$f$ must be 2-valent on the unbounded component
of $\mathbb{C} \backslash f(S)$
since real values in the unbounded component of
$\mathbb{C} \backslash f(S)$ are 2-valent.
It follows from that Lemma \ref{lem:n0preims} that
$f$ must be 2:1 on the unbounded component of
$\mathbb{C} \backslash f^{-1}(f(S))$.
\end{example}

It is nicest when our partitioning set
$f(S) \cup C(f)$ has empty interior.
Note that the results above are vacuous if 
the partitioning set fills the complex plane.
For example,
$f(z) = e^z$ is harmonic in $\mathbb{C}$.
It is easy to see that $C(f) = \mathbb{C}$.
In this case, $Val(f, w) = \infty$ if $w \neq 0$
and $Val(f, 0) = 0$.

\begin{remark}
If $z_0 \in R \subseteq \mathbb{C}$,
where $R$ is open,
then there exists a simply connected, open set 
$U \subseteq R$ such that $z_0 \in U$
(for example, 
$B(z_0, \delta)$ for $\delta > 0$ sufficiently small.)
Suppose that $f$ is harmonic in $R$.
Since $U$ is simply connected,
there exist functions $h$ and $g$ holomorphic
in $U$ such that
$f = h + \overline{g}$ in $U$.
We have $J_f = |h'|^2 - |g'|^2$ in $U$.
Thus,
$z \in S \cap U$ iff $|h'(z)| = |g'(z)|$.
If $R$ is a simply connected, open set,
then we can take $U = R$. 
\end{remark}

\begin{lemma}
\label{lem:nullpartifsnotempty}
Let $f$ be harmonic in a simply connected, open set 
$R \subseteq \mathbb{C}$.
If $int\ S \neq \varnothing$,
then $f(z) = \alpha\ +\ \beta u(z)$,
where $u$ is a real harmonic function 
and $\alpha$, $\beta$ are complex constants.
Moreover,
$S = R$ and $Val(f, w_0) = \infty$ for every
$w_0 \in f(R)$.
\end{lemma}
\begin{proof}
Since $R$ is simply connected,
there exist functions $g$ and $h$ holomorphic in $R$
such that
$f = h + \overline{g}$ in $R$.

First, suppose that the zeros of $g^\prime$ are
not isolated in $R$.
By the identity theorem, 
we must have $g^\prime \equiv 0$ in $R$;
thus $g$ is constant in our connected set $R$.
Since $int\ S \neq \varnothing$,
there exists a non-empty open disc 
$U \subseteq S$.
Since $g' \equiv 0$,
$h^\prime \equiv 0$ in $U$.
Since $h$ is constant in $U$, 
$h$ is constant in our connected set $R$.
Hence, $f = h + \overline{g}$ is constant
and maps $R$ to a single point;
that point thus has infinite valence.
The representation of $f$ follows trivially.
Also note that $S = R$.

Now, suppose that the zeros of $g^\prime$ are isolated.
We may then choose a non-empty open set 
$U \subseteq S$
such that $g^\prime(z) \neq 0$ for all $z \in U$.
Then $\psi(z) = h^\prime(z) / g^\prime(z)$ is
a holomorphic function on $U$ and $|\psi| \equiv 1$
on $U$.
Take a non-empty closed disc $K \subset U$.
By the maximum modulus principle, 
since $|\psi| \equiv 1$ on $K$,
$\psi \equiv \lambda$ on $K$ for some
unimodular constant $\lambda$.
Moreover,
$h^\prime - \lambda g^\prime \equiv 0$ on $int\ K$.
Thus, $h^\prime - \lambda g^\prime \equiv 0$ on 
$R$ and $h^\prime(z) = \lambda g^\prime(z)$
for all $z \in R$.
We thus have that $S = R$.
Further,
for some constant $\alpha \in \mathbb{C}$,
$h = \alpha + \lambda g$. 
Choose $\tau$ such that 
$\tau^2 = \lambda$.
Then $f = h + \overline{g} = 
\alpha + 2 \tau\ Re(\tau g) =
\alpha + \beta u$
in $R$,
where $\beta = 2 \tau$
and $u = Re(\tau g)$ is a real-valued harmonic
function in $R$.
Choose $w_0 \in f(R)$.
Then we may choose $z_0 \in R$ such that 
$f(z_0) = w_0$.
Thus, $\alpha = w_0 - \beta u(z_0)$.
Consider $v(z) = u(z) - u(z_0)$.
Since the zeros of a real-valued harmonic function
are not isolated,
$v$ has an infinite number of distinct
zeros in $R$;
hence $Val(f, w_0) = \infty$.
\end{proof}

\begin{remark}
The preceding lemma can also be found 
in~\cite{L:light}.
If $f$ is harmonic in a simply connected, open 
set in $\mathbb{C}$ and if $int\ S \neq \varnothing$,
then the preceding lemma implies that
$f(R)$ is either a point or a subset of a line.
Thus, $int(f(S) \cup C(f)) = \varnothing$.
Moreover, 
if $w \in \mathbb{C} \backslash (f(S) \cup C(f))$,
then $Val(f, w) = 0$.
\end{remark}

\begin{lemma}
\label{lem:snotempty}
Let $f$ be harmonic in an open set 
$R \subseteq \mathbb{C}$.
If $int\ S \neq \varnothing$,
then there exists a point 
$w_0 \in f(R)$
such that $Val(f, w_0) = \infty$.
Moreover,
there exists an open, non-empty set $R_1$, 
with $R_1 \subseteq S \subseteq R$,
such that
$f(R_1)$ is either a point or is a subset of a line.
\end{lemma}
\begin{proof}
Choose $z_0 \in int\ S$.
Then we can find $\delta > 0$
such that $B(z_0, \delta) \subseteq int\ S$.
Let $R_1 = B(z_0, \delta)$.
Then $R_1$ is a non-empty, simply connected, open set
and we may apply Lemma \ref{lem:nullpartifsnotempty}
to $f|_{R_1}$.
The result follows.
\end{proof}

\begin{corollary}
\label{cor:goodpart}
Let $f$ be harmonic in an open set 
$R \subseteq \mathbb{C}$.
If $Val(f, w)$ is finite for each $w \in \mathbb{C}$,
then $int(f(S) \cup C(f)) = \varnothing$.
\end{corollary}
\begin{proof}
Since $f$ is $C^\infty$ on $R$,
$f(S)$ has Lebesgue measure 0 by Sard's Theorem;
hence, 
$f(S)$ has empty interior.
By Lemma \ref{lem:snotempty},
since we assume that $Val(f,w)$ is finite for each
$w \in \mathbb{C}$,
we must have $int\ S = \varnothing$.
Since $S$ is relatively closed in $R$, 
we have that $S$ is nowhere dense in $R$.
Since we have assumed that
$\{w: Val(f,w) = \infty\} = \varnothing$,
$int (C(f)) = \varnothing$ by Theorem \ref{thm:c1dense}.
The result follows from Remark \ref{rem:emptyint} and 
Theorem \ref{thm:c1dense}.
\end{proof}

The results above relate to a theorem of 
T. J. Suffridge and J. W. Thompson~\cite{ST:partition}.
We need a definition before stating the theorem.

\begin{definition} 
[Suffridge and Thompson]
Assume $f$ is harmonic in 
$\{z: 0 < |z-\alpha| < r\}$ for some $r > 0$.
Define $\alpha$ to be a pole of $f$ provided
$lim_{z \rightarrow \alpha} |f(z)| = \infty$.
\end{definition}

\begin{theorem} 
[Suffridge and Thompson]
\label{thm:stpartition}
Let $f$ be non-constant and harmonic,
except for a finite number of poles,
in a simply connected region $D$ in $\mathbb{C}$.
Let $C$ be a Jordan curve contained within $D$ not 
passing through a pole.
Let $\Omega$ be the region bounded by $C$.
Let $X$ be the exceptional set of $f$ in $D$.
Then the image space is partitioned into open
components by $f(C \cup X)$ such that all values in
a given component are assumed the same number of times
in $\Omega$.
\end{theorem}

\noindent
The ``exceptional set of $f$'' consists of points 
in $D$
where the 
dilatation of $f$
(defined as $\overline{f_{\overline{z}}}/ f_z$)
is unimodular.
The proof of Theorem \ref{thm:stpartition} 
in~\cite{ST:partition} assumes that the dilatation
is analytic in $D$;
this occurs if $S$ has no isolated critical points.
When $S$ has no isolated critical points,
$S = X$.
As Example \ref{ex:isolcrit} below demonstrates,
we must include the images of isolated 
critical points in the partitioning set.

We can correct this difficulty and extend
Theorem  \ref{thm:stpartition} to the case of less
restrictive conditions on $D$ and its boundary
by using the partitioning results above, as follows:

\begin{corollary}
\label{cor:stgen}
Let $f$ be harmonic in a bounded, open set 
$R \subset \mathbb{C}$.
Suppose that $P = int (\overline{R}) \backslash R$
is empty or consists of a finite number of poles of 
$f$.
Suppose also that $f$ can be continuously
extended to $\partial R \backslash P$.
Then $f(S) \cup f(\partial R \backslash P)$
partitions $\mathbb{C}$ into regions of constant
valence.
\end{corollary}
\begin{proof}
If we have a sequence 
$\{z_n\} \subset R$ approaching $z \in P$,
then $f(z_n) \rightarrow \infty$.
Thus,
since $R$ is bounded,
$C(f) = \{w: \exists \{z_n\} \subset R$ such that
$f(z_n) \rightarrow w$ and $z_n \rightarrow
z_0 \in \partial R \backslash P\}$.
Since $f$ can be continuously extended to
$\partial R \backslash P$,
$C(f) = f(\partial R \backslash P)$.
Then apply Theorem \ref{thm:harmpart}.
\end{proof}

\begin{example}
\label{ex:isolcrit}
$f(z) = \frac{1}{3}\ z^3 - \frac{1}{2}\ \overline{z}^2$

In this example,
we will be looking at $f$ restricted to $R = B(0, 1)$.
This example illustrates why the image
($w = 0$)
of the isolated critical point 
($z = 0$)
must be included in the partitioning set.
We will look at $f$ as a function in 
$\mathbb{C}$ in Example \ref{ex:plane_isolcrit}
below.

Looking at $f$ in $B(0,1)$,
calculations show that $S = \{0\}$
and that $C(f) = f(\{z: |z| = 1\})$.
The partitioning set is graphed in 
Figure \ref{fig:isolcrit_imaglabel}
(see Example \ref{ex:plane_isolcrit} below)
where the cusps are identified by numbers.
Note that the partitioning set is the same set of 
points as in Example \ref{ex:plane_isolcrit}.

We will explicitly find the preimages of
$w = 0$ 
(Claim 2)
and show that the valence is not
constant in a sufficiently small neighborhood of
$w = 0$
(Claim 3.)
It is helpful to rewrite $f$ as a function
of $(x,y)$ where $z = x + iy$: 

\begin{center}
$f(x, y) = \frac{1}{3}x^3 - x y^2 - \frac{1}{2}x^2
+ \frac{1}{2}y^2 + i y\ (x^2 - \frac{1}{3}y^2 + x)$.
\end{center}

\noindent
\textbf{Claim 1:}
Let $y^2 = 3 x\ (x + 1)$.  
Then $(x, y) \in B(0,1)$ iff $0 \le x < 1/4$.
\begin{proof}[Outline of proof]
Note that $(x, y) \in R = B(0,1)$ 
iff $0 \le x^2 + y^2 < 1$.
Using our condition on $y^2$, 
we require $-1 \le 4 x^2 + 3 x - 1 < 0$.
We note that $\varphi(x) = 4 x^2 + 3 x - 1$ has
zeros at $x = -1$, $1/4$ and
attains its minimum at $x = -3/8$.
So we require $-1 < x < 1/4$.
Now,
$\varphi(0) = -1$ and $\varphi(-3/8) < -1$.
We also need $y^2 = 3 x\ (x + 1) \ge 0$ and
that forces us to require $x \ge 0$.
\end{proof}

\noindent
\textbf{Claim 2:}
$Val(f|_R, 0) = 1$.
\begin{proof}[Outline of proof]
To satisfy $Im\ f(x, y)\ =\ 0$, 
we must have $y = 0$ or $y^2 = 3 x\ (x + 1)$.
We must also satisfy $Re\ f = 0$.
We can rewrite this as 
$(3 - 6 x)\ y^2 + x^2\ (2 x -3) = 0$.
It is easy to check that this equation has no
solutions when $x = \frac{1}{2}$.
Thus,
we must have $y^2 = \{x^2\ (2 x - 3)\} / (6 x - 3)$.
To get $f = 0$ when $y = 0$,
we must have 
$x = 0$ or $\frac{3}{2}$.
To get $f = 0$ when  
$y^2 = 3 x\ (x + 1)$ ,
we must have
$x = \frac{3}{8}(-1 \pm \sqrt{5})$;
neither value of $x$ satisfies Claim 1.
Also note that the corresponding values of
$y^2$ for both choices of $x$ are positive,
so these correspond to preimages of $0$ outside
of $R$.
So, $w = 0$ has six distinct preimages in 
$\mathbb{C}$,  
but only one preimage $R = B(0,1)$;
namely, $z = 0$.
\end{proof}

\noindent
\textbf{Claim 3:}
Choose $|\epsilon| > 0$, sufficiently small.
Then $Val(f|_R, \epsilon) = 2$.
\begin{proof}[Outline of proof]
As in Claim 2 above,
$Im\ f = 0$ requires $y = 0$
or $y^2 = 3 x\ (x + 1)$.
We also require that $Re\ f = \epsilon$.

First, let $y = 0$.
Then $Re\  f = \epsilon$ becomes $\varphi(x) = 0$,
where 
$\varphi(x) = 
\frac{1}{3}x^3 - \frac{1}{2}x^2 - \epsilon$
and $-1 < x < 1$ 
(to get a preimage in $R$.)
Now $\varphi'(x) = 0$ at $x = 0, 1$;
$\varphi(0) = -\epsilon$ and 
$\varphi(1) = -(\epsilon + \frac{1}{6})$.
Note that $\varphi(-1) = -(\frac{5}{6} + \epsilon)$.
Also $\varphi''(x) = 0$ at $x = 1/2$.
An elementary calculus argument shows that
we have no solutions for $\varphi(x) = 0$ in $(-1, 1)$
if $0 < \epsilon < 1/6$.
Also,
if $-\frac{1}{6} < \epsilon < 0$,
$\varphi$ has two distinct zeros in $(-1, 1)$.
Thus, for $|\epsilon|$ sufficiently small:
\begin{itemize}
\item
If $\epsilon < 0$,
two distinct preimages of $w = \epsilon$ in
$R$ on the real axis.
\item
If $\epsilon > 0$,
no preimages of $w = \epsilon$ in $R$ on the real
axis.
\end{itemize}

Now, let $y^2 = 3 x\ (x + 1)$.
We may assume $y \ne 0$,
since that case has already been handled.
Then $Re\ f = \epsilon$ becomes $\varphi(x) = 0$,
where 
$\varphi(x) = 
\frac{8}{3}x^3 + 2 x^2 - \frac{3}{2} x +
\epsilon$.
By Claim 1,
we also require $0 \le x < 1/4$.
The zeros of $\varphi'(x)$ are $x = 1/4, -3/4$,
while $\varphi''(x) = 0$ at $x = -1/4$.
Notice that $\varphi(0) = \epsilon$,
$\varphi(1) = \frac{19}{6} + \epsilon > 0$
for $|\epsilon|$ sufficiently small,
and $\varphi(1/4) = -\frac{5}{24} + \epsilon < 0$
for $|\epsilon|$ sufficiently small.
Thus, for $|\epsilon|$ sufficiently small we
have:
\begin{itemize}
\item
If $\epsilon < 0$, 
$\varphi(x)$ has no zeros in $[0, 1/4)$.
So, no preimages of $w = \epsilon$ in $R$ off the
real axis.
\item
If $\epsilon > 0$,
$\varphi(x)$ has one zero in $[0, 1/4)$;
that root is not zero.
Since  $y^2 = 3 x\ (x + 1) \ne 0$ for
$x \in (0, 1/4)$,
two distinct preimages of $w = \epsilon$ in $R$ off 
the real axis.
\end{itemize}

\noindent
The claim follows by combining the two cases.
\end{proof}
\end{example}

\noindent
When does the partitioning set in the preceding
corollary have empty interior?

\begin{lemma}
Let $f$ be harmonic in a bounded open set 
$R \subset \mathbb{C}$.
Suppose that $P = int (\overline{R}) \backslash R$
is empty or consists of a finite number of poles
of $f$.
Suppose also that $f$ has a $C^1$ extension to
some open set $R_1$,
where $R_1 \supset \overline{R} \backslash P$.
Then $f(S) \cup f(\partial R \backslash P)$ 
has empty interior.
\end{lemma}
\begin{proof}
The following proof was suggested by D. Sarason.
Let $F$ denote our $C^1$ extension of $f$ to $R_1$
and $S_F$ denote the critical set of $F$.
Let
\begin{eqnarray*}
A &=& \partial R \backslash (S_F \cup P) \\
B &=& S \cup (S_F \cap (\partial R \backslash P))
\end{eqnarray*}
Then $S \cup (\partial R \backslash P) = A \cup B$.
If $F(A) \cup F(B)$ is a countable union
of nowhere dense sets,
then the result follows from the Baire Category
Theorem.

We claim that $F(B)$ is a countable union of
nowhere dense sets.
By construction, 
$B \subseteq S_F$;
hence $F(B) \subseteq F(S_F)$ and it is enough
to show that $F(S_F)$ is a countable union of
nowhere dense sets.
$F(S_F)$ has measure 0 by Sard's Theorem.
Since $S_F$ is closed in $R_1$,
$S_F$ is a countable union of compact sets.
Therefore $F(S_F)$ is a countable union of compact
sets,
each of which has measure 0,
hence each nowhere dense.

We also claim that
$F(A)$ is a countable union of nowhere dense sets.
We note that $F$ is a local homeomorphism at each
$z \in A$,
so for each $z \in A$, there exists $\epsilon > 0$
such that $F|_{B(z, 2 \epsilon)}$ is a homeomorphism.
In particular,
$B(z, 2 \epsilon) \cap (S_F \cup P) = \varnothing$.
Then $\mathcal{C} = \{B(z, \epsilon): z \in A\}$
is an open cover of $A$.
By Lindel\"{o}f's Theorem,
$\mathcal{C}$ has a countable subcover,
say $\{B(z_j, \epsilon_j): z_j \in A\}$.
By construction,
$K_j = \overline{B(z_j, \epsilon_j) \cap A}
\subseteq A$
and we see that $A = \cup K_j$.
We note that $K_j \subseteq \partial R$ is compact and 
has empty interior.
Since $F$ is a homeomorphism on $K_j$,
$F(K_j)$ is nowhere dense.
The claim follows by noting that 
$F(A) = \cup F(K_j)$.
\end{proof}

\section {Application to 
harmonic mappings in the entire complex plane}
\label{sec:entire}
If $f$ is harmonic in the entire complex plane,
$C(f) = C(f, \infty)$.
Also, 
we can write $f$ as $f = h + \overline{g}$
where $h$ and $g$ are entire (holomorphic) functions.
This allows us to say more about the behavior
of $f$ when $S$ has non-empty interior.
In particular,
we can apply 
Lemma \ref{lem:nullpartifsnotempty} as follows:

\begin{lemma}
\label{lem:plasnotempty}
Let $f$ be harmonic on $\mathbb{C}$.
If $int\ S \neq \varnothing$, 
then $f(\mathbb{C})$ is either a single point or
a line.
In either case,
$S = \mathbb{C}$ and there exists a point 
$w_0 \in \mathbb{C}$
such that $Val(f, w_0) = \infty$.
\end{lemma}
\begin{proof}
Let $f = h + \overline{g}$, where $g$ and $h$ are
entire (holomorphic) functions.
The conclusions about $S$ and the existence
of a point of infinite valence follow directly
from Lemma \ref{lem:nullpartifsnotempty}.

Now we return to the proof of 
Lemma \ref{lem:nullpartifsnotempty}.
If the zeros of $g'$ are not isolated,
then $f(\mathbb{C})$ is a single point.
If the zeros of $g'$ are isolated,
$g$ is a non-constant entire function and
there exists a unimodular constant $\lambda$
such that $h' \equiv \lambda g'$.
Choosing $\tau$ such that $\tau^2 = \lambda$
and choosing $z_0$ such that $f(z_0) = w_0$,
we have 
$f = h + \overline{g} = \alpha + 2 \tau\ Re(\tau g)$
where $\alpha = w_0 - 2 \tau\ Re(\tau g(z_0))$.

Since $\tau$ is non-zero and $g$ is a non-constant
entire function,
$\varphi(z) = Re(\tau g(z))$ is a non-constant, 
real harmonic 
function in the entire plane.
By Liouville's Theorem,
$\varphi(z)$ can be bounded neither above nor below.
Noting that
the range of $\varphi$ is connected,
we see that the range of $\varphi$ is the entire real
axis.
Hence $f(\mathbb{C})$ is a line when the zeros
of $g'$ are isolated.
\end{proof}

\begin{lemma}
\label{lem:planullpart}
Let $f$ be harmonic on $\mathbb{C}$.
Suppose that $int\ S \neq \varnothing$.
Then $int(f(S) \cup C(f,\infty)) = \varnothing$.
Moreover, if $w \in \mathbb{C} \backslash
(f(S) \cup C(f,\infty))$, 
then $Val(f, w) = 0$.
\end{lemma}
\begin{proof}
By Lemma \ref{lem:plasnotempty}, 
$S = \mathbb{C}$ and $f(\mathbb{C})$
is either a point or a line
(each is a closed set with empty interior 
in $\mathbb{C}$.)
Clearly, $C(f,\infty) \subseteq f(S) = f(\mathbb{C})$.
Hence, $f(S) \cup C(f,\infty)$ has empty interior.
Also, if $w \in \mathbb{C} \backslash
(f(S) \cup C(f,\infty))$,
then $w \in \mathbb{C} \backslash f(\mathbb{C})$
and thus $Val(f,w) = 0$. 
\end{proof}

We also can say more about when the partitioning
set $f(S) \cup C(f, \infty)$ has empty interior.
We first recall a result of Wilmshurst~\cite{W:thesis}:

\begin{theorem}[Wilmshurst]
\label{thm:wfinval}
Let $f(z) = \overline{g(z)} + h(z)$ be a function 
harmonic in the (entire) complex plane.
If $ \lim_{z \rightarrow \infty} 
f(z) = \infty$,
then $f$ has finitely many zeros.
\end{theorem}

Note that $ \lim_{z \rightarrow \infty} 
f(z) = \infty$
implies that $Val(f,w)$ is finite for each $w$
(note that this does not imply that $f$ has finite
valence.)
We also note that if $ \lim_{z \rightarrow \infty} 
f(z) = \infty$,
then $C(f, \infty) = \varnothing$.
Since $f(S)$ has empty interior by Sard's Theorem,
we have that $f(S) \cup C(f, \infty) = f(S)$ has empty
interior.
In particular,

\begin{corollary}
\label{cor:inflimpart}
If $f$ is a harmonic mapping on $\mathbb{C}$ such that
$ \lim_{z \rightarrow \infty} 
f(z) = \infty $,
then $ f(S) $ partitions $ \mathbb{C} $ 
into non-empty regions of constant valence.
\end{corollary}

\noindent
Corollary \ref{cor:inflimpart} obviously applies
if $f(z) = p(z) - \overline{q(z)}$,
where $p$ and $q$ are polynomials in $z$
of different degree.
D. Bshouty, W. Hengartner and 
M. Na\-ghi\-bi-Bei\-dokh\-ti
in\-~\cite{BHN:dilat} use approximation theory
to show that there exist harmonic mappings on 
$\mathbb{C}$ which are not polynomials such that 
$ \lim_{z 
\rightarrow \infty} f(z) = \infty$.
An explicit example of such a function is due
to D. Sarason;
consider $f(z) = Re(e^{-z^2}) + z$.
Calculations also show that $f$ is onto $\mathbb{C}$.
While $Val(f, w)$ is finite for each $w$,
$Val(f)$ is infinite.

We will show that $f(S) \cup C(f, \infty)$ has
empty interior when $f$ is a harmonic polynomial,
even when the Corollary \ref{cor:inflimpart}
does not apply.
To do so, 
we need to state several results.
We start with a result one can find 
in~\cite{W:paper}:

\begin{theorem}[Wilmshurst]
\label{thm:wconstarc}
Let a function $f$ be harmonic in some domain $D$ 
and have a 
sequence of distinct zeros $\{z_m\}$ converging
to some point $z^*$ in $D$.
Then $f(z) \equiv 0$ on some simple analytic arc 
containing $z^*$ as an interior point. 
Further,
there are at most finitely many such arcs unless
$f(z) \equiv 0$ in $D$.
\end{theorem}

We also note two forms of B\'{e}zout's Theorem
found in Wilmshurst ~\cite{W:paper}:

\begin{theorem}[Wilmshurst: B\'{e}zout's Theorem]
\label{thm:bezout}
Let $A$ and $B$ be relatively prime polynomials 
in the real
variables $x$ and $y$ with real coefficients,
and let $deg\ A = n$ and $deg\ B = m$.
Then the two algebraic curves $A(x, y) = 0$ and
$B(x,y) = 0$ have at most $mn$ points in common.
\end{theorem}

\begin{theorem}[Wilmshurst: 
Alternate form of B\'{e}zout's Theorem]
\label{thm:waltbezout}
Let $A$ and $B$ be polynomials in the real
variables $x$ and $y$ with real coefficients.
If $deg\ A = n$ and $deg\ B = m$, 
then either $A$ and $B$ have at most $mn$ common zeros
or have infinitely many common zeros.
\end{theorem}

\begin{remark}
\label{rem:algcurve}
B\'{e}zout's Theorem also provides information
about a real algebraic curve in the plane.
Let $P(x, y)$ be a polynomial in the real variables
$x$ and $y$ with real coefficients.
There exist irreducible polynomials in $x$ and $y$
with real coefficients, say 
$A_1(x, y), \ldots, A_k(x, y)$,
such that $P = A_1^{n_1} \ldots A_k^{n_k}$.
Let $V = \{(x, y): P(x, y) = 0\}$. 
Then $V = \cup_{j = 1}^{k} V_j$,
where $V_j = \{(x, y): A_j(x, y) = 0\}$.
Let $D_j = \{(x, y): \partial A_j / \partial y = 0\}$.
Since the $A_j$ are relatively prime,
if $i \neq j$,
then $V_i \cap V_j$ is finite by 
Theorem \ref{thm:bezout}.
If $\partial A_j / \partial y \equiv 0$,
then $V_j$ consists of a finite number
(possibly zero)
of vertical lines.
If $\partial A_j / \partial y$ is not identically 
zero,
then $V_j \cap D_j$ is finite
by Theorem \ref{thm:bezout} 
since $A_j$ is irreducible.
In that case,
we can divide the plane into a finite number of
vertical strips
such that $\partial A_j / \partial y \neq 0$
for points in $V_j$ in the interior of 
a given vertical strip.
By the implicit function theorem,
using analytic continuation,
we see that $V_j$ consists of a finite number
(possibly zero)
of non-intersecting curves in the interior
of a given 
vertical strip\cite[p. 249-253]{S:primer} and
possibly a finite number of vertical lines.
If $V_j$ contains a closed loop,
we see that $V_j$ must contain a simple loop.
A closed loop in $V$ is either a closed loop
in one of the $V_j$ or is formed by joining
arcs in $V$ terminating at points
in $V_i \cap V_j$ where $i \neq j$.
Since there are finitely many points in 
the $V_i \cap V_j$,
we see that if $V$ contains a closed loop,
then $V$ contains a simple loop.
\end{remark}

We also need a fact about harmonic polynomials 
from~\cite{BC:brelotchoquetlem}:

\begin{lemma} [M. Brelot and G. Choquet]
\label{lem:brelotchoquet}
If a real harmonic polynomial $p$ in two variables 
has a nonconstant factor $q$, 
then the zeros of $q$ are not isolated.
\end{lemma}

\begin{lemma}
\label{lem:constpolyonloop}
Let $f$ be a harmonic polynomial.
Let $\gamma$ be a closed loop in $S$.
If $S$ has empty interior,
then $f$ cannot be constant on $\gamma$.
\end{lemma}
\begin{proof}
Since $S$ has empty interior, 
we can view $S$ as a real algebraic curve.
By Remark \ref{rem:algcurve},
if $S$ contains a closed loop $\gamma$,
we can find a simple loop in $\gamma$.
Without loss of generality,
we may suppose that $\gamma$ is a Jordan curve.
Then $\gamma$ encloses some bounded, 
connected region, say $R$,
by the Jordan curve theorem.

Assume that $f$ is constant on $\gamma$.
Let $u = Re\ f$ and $v = Im\ f$.
Then $u$ and $v$ are real harmonic functions
in $\mathbb{C}$.
Since $f$ is constant on $\gamma = \partial R$,
so are $u$ and $v$.
Hence, $u$ must be constant in $R$
by the minimum and maximum principles for
real harmonic functions.
Similarly,
$v$ is constant in $R$.
Thus,
$f$ is constant in $R$ and hence in $\mathbb{C}$. 
We have $S = \mathbb{C}$,
a contradiction.
\end{proof}

\noindent
The preceding proof uses an argument similar
to that used by Wilmshurst
(pages 42-43 in~\cite{W:thesis})
to show that the zeros of a harmonic polynomial
are isolated under certain conditions.

\begin{corollary}
\label{cor:partpoly}
If $f$ is a harmonic polynomial,
then $int (f(S) \cup C(f, \infty)) = \varnothing$.
\end{corollary}
\begin{proof}
The result is obvious if $f$ is constant,
so we will assume that $f$ is a non-constant
harmonic polynomial.
Following Wilmshurst~\cite{W:paper},
finding all $z$ such that $f(z) - w = 0$ is
equivalent to finding the common zeros of
$Re(f-w)$ and $Im (f- w)$ in $\mathbb{R}^2$,
identifying $\mathbb{C}$ with $\mathbb{R}^2$ in
the obvious way.

Let $n = Max\{deg\,(Re\,f), deg\,(Im\,f)\}$.
By Theorem \ref{thm:waltbezout},
$Re (f - w)$ and $Im (f - w)$ have either at most
$n^2$ common zeros or have infinitely many 
common zeros.
Hence $Val(f) \le n^2$ or there exists
$w_0 \in \mathbb{C}$ such that $Val(f, w_0) = \infty$.

If $Val(f)$ is finite, 
the result follows from Corollary \ref{cor:goodpart}.
So we may suppose that there exists
$w_0 \in \mathbb{C}$ such that $Val(f, w_0) = \infty$.
If $S$ has non-empty interior,
we are in the simple case where
$f$ maps $\mathbb{C}$ to a point or a line by
Lemma \ref{lem:plasnotempty},
and the result follows.
So, 
suppose also that $S$ has empty interior.
Since $S$ is closed in $\mathbb{C}$,
we also have that $S$ is nowhere dense.
Since $f$ is harmonic,
$f(S)$ has empty interior by Sard's Theorem.
We want to show that
$W_\infty = \{w: Val(f, w) = \infty\}$ 
is nowhere dense.
If that is true,
the result follows 
from Theorem \ref{thm:c1dense} 
and Remark \ref{rem:emptyint}.
We will show that $W_\infty$ is finite.

Since $Val(f, w_0) = \infty$,
there are infinitely many distinct values of $z$
such that $f(z) = w_0$.
By B\'{e}zout's Theorem (Theorem \ref{thm:bezout}),
$Re (f - w_0)$ and $Im (f - w_0)$ have a non-trivial
common factor.
Choose $z_0$ such that $f(z_0) = w_0$.
By the Brelot-Choquet Lemma 
(Lemma \ref{lem:brelotchoquet}),
$z_0$ is not an isolated zero of $f(z) - w_0$.
Hence,
by Theorem \ref{thm:wconstarc},
there exists a simple analytic arc $\gamma$ such
that $f\mid_\gamma\ \equiv\ w_0$ and such that 
$z_0$ is in the interior of $\gamma$.
Since $f$ is not locally 1-1 at any point in the
interior of $\gamma$,
$int\ \gamma \subseteq S$.
Further,
we can apply Theorem \ref{thm:wconstarc}
to each of the endpoints of $\gamma$
to see that $\gamma \subseteq S$. 
If $\gamma$ contains a closed loop,
then $S$ has non-empty interior
by Lemma \ref{lem:constpolyonloop},
a contradiction.
So we may suppose that $\gamma$ does not contain
a closed loop.
Thus,
for each $w \in W_\infty$,
we can find 
a simple curve $\gamma$
in $S$ such that $f|_\gamma \equiv w$.

Since $f$ is a harmonic polynomial, 
$J_f(z)$ can be viewed as a polynomial in
$\mathbb{R}^2$ and Remark \ref{rem:algcurve}
applies to $S$.
Moreover,
if we extend $\gamma$ by repeated application
of Theorem \ref{thm:wconstarc},
we see that $\gamma$ extends to an unbounded curve
such that $f|_\gamma \equiv w_0$.
Otherwise,
if our extensions of $\gamma$ have a finite
limit point $z_0$,
by continuity,
$f(z_0) = w_0$ and we can extend $\gamma$ through
$z_0$,
a contradiction.
From Remark \ref{rem:algcurve},
we see that $\gamma$ can be extended to
a simple, unbounded curve in $S$.
From Remark \ref{rem:algcurve}, 
we see that $S$ contains at most finitely many 
simple, unbounded curves.
Hence $W_\infty$ contains a finite number of points and
the result follows.
\end{proof}

\noindent
Example \ref{ex:flatpoly} is an example
of a harmonic polynomial with $int\ S = \varnothing$ 
and a point with infinite valence.
In this example,
the preimages of this point are non-isolated
and lie on a curve (the line $Re\ z = Im\ z$.)

The situation is not so nice when $f$ is 
not a harmonic polynomial.
We will restrict our attention to the case
$R = \mathbb{C}$.
For example,
if $f$ is an entire (holomorphic) transcendental 
function,
$C(f, \infty)$ is the complex plane
by Picard's Theorem.
The monograph of M. Balk~\cite{Balk:book}
concerns the theory of polyanalytic and polyentire
functions;
results similar to Picard's Theorem are given
for polyentire functions
(see pages 107-109.)

Recall that $\varphi(z)$ is a polyanalytic function 
of order $n$
if it can be written as
$\varphi(z) = \Sigma_{k = 0}^{n-1}\  
a_k(z) \overline{z}^k$
where each $a_k(z)$ is a holomorphic function of $z$.
If each $a_k(z)$ is entire,
then $\varphi$ is polyentire.
Let $h$ and $g$ be holomorphic in some region $R$
and let 
$f = h + \overline{g}$.
Then the harmonic function $f$ can be thought of as 
a polyanalytic function 
(of countable order if $g$ is transcendental.)
If $g(z) = \sum_{k=0}^m b_k z^k$ 
is a non-constant polynomial,
we can think of $f$ as a polyanalytic function
of finite order by putting 
$a_0(z) = h(z) + \overline{b_0}$ and 
$a_k(z) \equiv \overline{b_k}$ 
for $1 \le k \le deg\ g$.

Let $P(z, \overline{z})$ denote a polynomial
in $z$ and $\overline{z}$.
Note that $P$ can include terms of the form
$z^n \overline{z}^m$.
Consider the following result from~\cite{Balk:book}
about polyentire functions:

\begin{theorem}[Balk]
\label{thm:polyentire}
Let $E(z)$ be an entire transcendental analytic
function,
and let $P(z, \overline{z})$ be an arbitrary
polyanalytic but not analytic polynomial.
Then the polyentire function
$F(z) = E(z) + P(z, \overline{z})$ assumes in the
complex plane $C$ every complex value $A$
(without any exceptions!)
The set $M(F; A)$ of all A-points of the
polyentire function $F(z)$ is for every $A$
unbounded and discrete in $C$.
\end{theorem}

If $f = h + \overline{g}$ is harmonic in $\mathbb{C}$,
the preceding result tells us that 
$C(f, \infty) = \mathbb{C}$ if $h$ is
transcendental and $g$ a non-constant polynomial.
Picard's Theorem takes care of the case when
$h$ is transcendental and $g$ a constant.
This gives us a necessary condition for
$f(S) \cup C(f, \infty)$ to have empty interior:

\begin{corollary}
\label{cor:badpart}
Let $f = h + \overline{g}$ be harmonic in 
$\mathbb{C}$, 
where $h$ and $g$ are entire.
If $f(S) \cup C(f, \infty)$ has empty interior,
then one of the following holds:
(i) $h$ and $g$ are both polynomials
in $z$
or
(ii) $h$ and $g$ are both entire transcendental
functions in $z$.
\end{corollary}

\noindent
Does the converse hold?
Corollary \ref{cor:partpoly} shows that the 
converse holds for harmonic polynomials.
We have seen some examples where $h$ and $g$
are both transcendental and $f(S) \cup C(f, \infty)$
has empty interior
(see  Example \ref{ex:transharm} and
the comments after Corollary \ref{cor:inflimpart}.)
Both examples are of the form 
$f(z) = z + Re\ h(z)$
where $h$ is an entire transcendental function.

\subsection{Joining regions for $f$ harmonic}
Let $f$ be a light harmonic function in $\mathbb{C}$.
Suppose that $R_1$ and $R_2$ are two distinct 
components of 
$\mathbb{C} \backslash 
f^{-1}(f(S) \cup C(f, \infty))$.
Suppose that that $R_1$ and $R_2$ 
share a common boundary arc.
Let $f(R_1) = \Omega_1$ and $f(R_2) = \Omega_2$.
If $\Omega_1$ and $\Omega_2$ are
each simply connected,
then $f$ is univalent in $R_1$ and in $R_2$
by Theorem \ref{thm:preimcover}.
Can we join $R_1$ and $R_2$ along the interior of
the shared
boundary arc to get a new region where
$f$ is univalent?

If the regions we want to combine all lie
in a convex domain,
a partial answer is given in~\cite{W:thesis}:

\begin{theorem}[Sheil-Small]
\label{thm:univconvex}
If $g(z)$ is an analytic function in the convex
domain $D$ and $|g'(z)| < 1$ in $D$ then
$\overline{z} + g(z)$ is univalent on $D$.
\end{theorem}

\noindent
Let $h(z)$ be holomorphic in some convex
domain $D$ and suppose that 
$|h'(z)| < 1$ in $D$.
It is easy to see that this result also
holds for $f(z) = h(z) + \lambda \overline{z}$
where $\lambda$ is a unimodular constant 
(apply the theorem to 
$f_1(z) = \overline{\lambda} f(z)$.)

We can apply this result to each of the bounded
components of $\mathbb{C} \backslash S$ in
Example \ref{ex:cubictwos} below to show
that $f$ is univalent on each of these two
bounded components.
However,
this result won't help us in the case that
a bounded component of $\mathbb{C} \backslash S$
is not convex;
consider for example the critical set in 
Example \ref{ex:dumbbell} below.
Also,
this result will not help us find univalent regions
in the unbounded component of 
$\mathbb{C} \backslash S$
nor will it help us with more complicated 
harmonic functions.

What do the results in the preceding sections
tell us about finding univalent regions?
By identifying $\mathbb{C}$ with $\mathbb{R}^2$
in the obvious way,
we may apply the preceding results.

\begin{lemma}
\label{lem:puncturepoints}
Let $f$ be a light harmonic function in $\mathbb{C}$.
Let $R$ be a bounded component of
$\mathbb{C} \backslash f^{-1}(f(S) \cup C(f, \infty))$.
Choose a component $\Omega \subseteq
\mathbb{C} \backslash (f(S) \cup C(f, \infty))$
such that $f(R) = \Omega$.
Suppose that $w_0$ is an isolated point in 
$\partial \Omega$.
If $f^{-1}(w_0)\ \cap\ \partial R\ \cap\ S = \varnothing$,
then $w_0$ has exactly $Val(f|_R, \Omega)$
distinct preimages in $\partial R$.
These preimages are isolated points in $\partial R$.
\end{lemma}
\begin{proof}
Note that $\Omega$ exists by Theorem 
\ref{thm:preimonto}.
We first check that 
$\zeta \in f^{-1}(w_0) \cap \partial R$ is isolated
in $\partial R$.
Suppose not.
We first note that $\zeta$ is isolated in 
$f^{-1}(w_0)$.
Otherwise, 
by Theorem \ref{thm:wconstarc},
there exists a non-degenerate arc $\gamma$
such that $f$ is constant on $\gamma$.
This contradicts $f$ being light;
hence $\zeta$ is isolated in $f^{-1}(w_0)$.
By Lemma \ref{lem:harmbdyim},
$f(\partial R) = \partial \Omega$.
Hence if $\zeta$ is not isolated in $\partial R$,
then there exists $\{\zeta_n\} \subseteq 
\partial R \backslash \{\zeta\}$
such that $\zeta_n \rightarrow \zeta$.
Since $\zeta$ is isolated in $f^{-1}(w_0)$,
$f(\zeta_n) \neq w_0$ for $n$ sufficiently large.
Since $\{f(\zeta_n)\} \subseteq \partial \Omega$
and $f(\zeta_n) \rightarrow w_0$ by continuity,
$w_0$ is not isolated in $\partial \Omega$,
a contradiction. 

We need to check that the conditions in 
Theorem \ref{thm:harmbdypreim} are satisfied.
Since $R$ is bounded,
$C(f|_R, \infty) = \varnothing$ and (i) holds.
By assumption,
(ii) holds.
Since $R$ is bounded and $f$ is harmonic
in $\mathbb{C}$,
we can look at a sequence in $\Omega$ converging
to $w_0$ and look at the corresponding sequence
of preimages in $R$ to see that
(iii) holds.
Since the preimage points of $w_0$ are isolated
in $\partial R$,
(iv) also holds.
The claim concerning the number of preimages  
follows from Theorem \ref{thm:harmbdypreim}.
\end{proof}

\begin{corollary}
\label{cor:fillpuncture}
Let $f$ be a light harmonic function in $\mathbb{C}$.
Let $\Omega$ be a component of 
$\ \mathbb{C} \backslash (f(S) \cup C(f, \infty))$.
Suppose that $(int\ \overline{\Omega}) \backslash
\Omega$ consists of a finite number of points.
Suppose that there exists a bounded component,
say $R$,
of $\ \mathbb{C} \backslash f^{-1}
(f(S) \cup C(f, \infty))$ 
such that $f(R) = \Omega$.
Suppose that every point in 
$(int\ \overline{R}) \backslash R$ is mapped to
$(int\ \overline{\Omega}) \backslash \Omega$.
If $(int\ \overline{R}) \backslash R$ contains no
point in $S$, 
then $f$ is $Val(f|_R, \Omega): 1$ in 
$int\ \overline{R}$.
Moreover,
if $int\ \overline{\Omega}$ is simply connected,
then $f$ is univalent in $int\ \overline{R}$.
\end{corollary}
\begin{proof}
Since $R$ is bounded,
we can apply Lemma \ref{lem:puncturepoints}
to $R$ to get that each of the isolated points
in $\partial \Omega$ has $Val(f|_R, \Omega)$
preimages in $\partial R$,
which are puncture points in $R$.
Thus,
each point in $int\ \overline{\Omega}$ has
exactly $Val(f|_R, \Omega)$ distinct
preimages in $int\ \overline{R}$,
the first claim.
Since $S \cap int\ \overline{R} = \varnothing$,
we see that $int\ \overline{R}$ is a 
$Val(f|_R, \Omega)$-fold covering of
$int\ \overline{\Omega}$
(see the proof of Theorem \ref{thm:preimcover}.)
The second claim follows from 
Theorem \ref{thm:munch8th4.5}.
\end{proof}

\begin{figure}
\includegraphics[scale=.95]{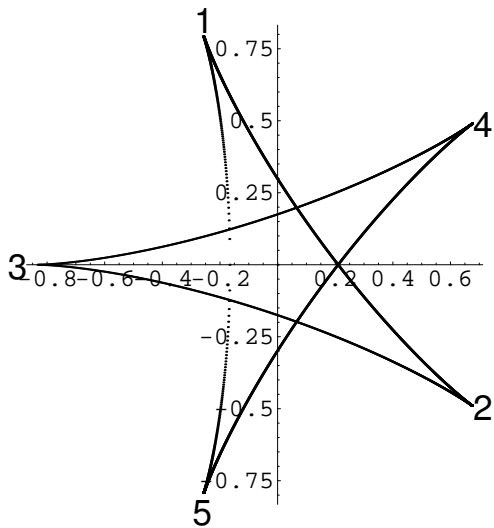}
\caption{Image of critical set for $f(z) = 
\frac{1}{3}\ z^3 - \frac{1}{2}\ \overline{z}^2$,
with labels.
Note that $0 \in f(S)$.}
\label{fig:isolcrit_imaglabel}
\end{figure}

\begin{figure}
\includegraphics[scale=1.25]{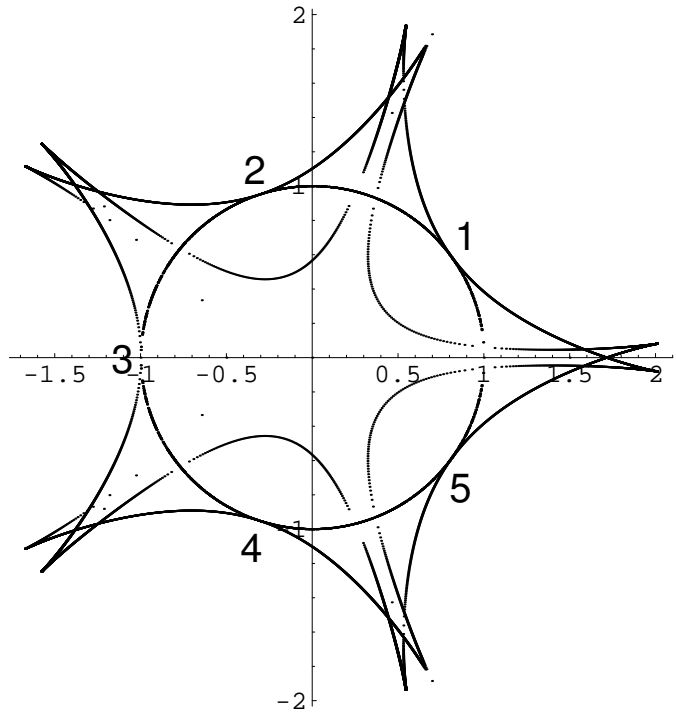}
\caption{Preimage of $f(S)$ for $f(z) = 
\frac{1}{3}\ z^3 - \frac{1}{2}\ \overline{z}^2$}
\label{fig:isolcrit_preim}
\end{figure}

\begin{example}
\label{ex:plane_isolcrit}
$f(z) = \frac{1}{3}\ z^3 - \frac{1}{2}\ \overline{z}^2$

We looked at this function
restricted to $\{z: |z| < 1\}$ in 
Example \ref{ex:isolcrit} above.
We will now look at $f$ as a function in $\mathbb{C}$.
Calculations show that $S = \{z: |z| = 1\} \cup \{0\}$
and that $C(f, \infty) = \varnothing$.
The partitioning set for $f$ will be $f(S)$,
which is graphed in 
Figure \ref{fig:isolcrit_imaglabel}.
We can select a few points in each region of
the partition and ask Mathematica to find the
preimages.
$Val(f, w) = 3$ in the unbounded
component of $\mathbb{C} \backslash f(S)$.
In each ``point of the star'' region,
$Val(f, w) = 5$.
$Val(f, w) = 7$ in the ``center of the star'' region
(recall that the origin is in $f(S)$ and is not in
this region.)
We have numbered the cusps in $f(S)$
in Figure \ref{fig:isolcrit_imaglabel}.
The critical points mapped to a given cusp
are labeled with the same number
in Figure \ref{fig:isolcrit_preim}.

Recall that the origin is an isolated critical point
of $f$.
From the calculations in Example \ref{ex:isolcrit},
we see that the origin has six distinct preimages.
From Figure \ref{fig:isolcrit_preim},
we see that the preimages other than the origin
lie in regions in the arms of the ``star''
(one in each such region.)
From the calculation in Example \ref{ex:isolcrit},
we see that $f$ is 2:1 in the component with
the origin as a puncture point.
Hence $f$ must be univalent in each of the five other
components which have $f^{-1}(0)$ as a puncture point.
It's not completely obvious from the figure
that the preimage of the origin is an isolated
point in each of the regions,
because the Mathematica routines plot single points.
However,
in a region where $f^{-1}(0)$ is clearly a puncture
point
(for example,
the arm in the lower right quadrant),
this serves as an example of Corollary
\ref{cor:fillpuncture}.
We can fill in the puncture point not in $S$
in each such component of the preimage to get
a simply connected region where $f$ is univalent.
\end{example}

\begin{theorem}
\label{thm:jointwo}
Let $f$ be a light, harmonic function in $\mathbb{C}$.
Let $R_1$ and $R_2$ be disjoint bounded
components of 
$\ \mathbb{C} \backslash 
f^{-1}(f(S) \cup C(f, \infty))$.
Suppose that $R_1$ and $R_2$ each have at most a
finite number of puncture points.
Suppose that for each component $\Omega \subseteq
\mathbb{C} \backslash (f(S) \cup C(f, \infty))$,
$\partial \Omega \cap int\ \overline{\Omega}$ 
consists of a finite number of points.
Suppose that $f$ is univalent on $R_1$ and on $R_2$.
Let $\gamma$ be a non-degenerate, simple arc
in $\partial R_1 \cap \partial R_2$.
Suppose that $f(int\ \gamma)$ is a simple arc.
If $S\ \cap\ int\ \gamma = \varnothing$,
then
\begin{enumerate}
\item
$f(R_1)$ and $f(R_2)$ are distinct components
of $\ \mathbb{C} \backslash (f(S) \cup C(f,\infty))$.
\item
$f$ is a univalent on 
$R_1\ \cup\ R_2\ \cup\ int\ \gamma$. 
\end{enumerate}
\end{theorem}
\begin{proof}
This follows from Theorem \ref{thm:jointwoc1}
by associating $\mathbb{R}^2$ and $\mathbb{C}$
in the usual way.
\end{proof}

\begin{remark}
\label{rem:sharedarc}
Note that if $int\ \gamma$ contains a non-degenerate
subarc in $S$,
then there exists $z_0 \in int\ \gamma$ such
that $f_{z_0} \sim z, \overline{z}$
(see the proof of Theorem \ref{thm:lyzval}.)
Thus,
$f(R_1) \cap f(R_2) \neq \varnothing$.
By Theorem \ref{thm:preimonto},
$f(R_1) = f(R_2)$.
\end{remark} 

Extending these ideas to joining more than
two adjacent regions is more difficult.
Applying these results requires assumptions
as to whether pairs of components 
(both in the domain and in the range)
have boundaries intersecting in more than one
simple arc.
Even if we join two regions, 
it could be the case that the images of the
two regions meet in more than one simple arc;
joining along more than one of these arcs could
cause trouble.

If we have components as in Theorem \ref{thm:jointwo},
we can fill in puncture points if the conditions
in Corollary \ref{cor:fillpuncture} are satisfied
to get a simply connected region of univalence. 
The following theorem of M. Ortel and 
W. Smith~\cite{OrtelSmith:univalence}
will then help us join more than two univalent
components,
provided that the combined components in both
the domain and range
are simply connected.

\begin{theorem}[Ortel and Smith]
\label{thm:OSlcluniv}
If $\Omega_1$ and $\Omega_2$ are open,
connected, and simply connected subsets of the
complex plane, 
$f: \Omega_1 \rightarrow \Omega_2$ is continuous,
surjective, and locally univalent,
$N \in \{1,2,3, ...\}$ and 
$\# f^{-1}(w) \in \{1, N\}$ for all $w \in \Omega_2$,
then $f$ is univalent.
\end{theorem}

For example,
suppose that we have three disjoint 
simply connected components of
$\mathbb{C} \backslash f^{-1}(f(S) \cup C(f, \infty))$ 
where $f$ is univalent and 
the image of each component is simply connected.
Suppose also that the conditions of
Theorem \ref{thm:jointwo} hold pairwise.
We thus have either two or three distinct image 
components.
$f$ is univalent in the interior of each of the 
two boundary arcs used for joining regions.
By construction,
at least one of the image regions has
preimages in only one of the original components.
Thus, 
a point in the image of the combined image region must
have either one or two preimages in 
the combined region.
Univalence in the combined region follows
from the theorem of Ortel and Smith,
provided that we show that the combined
image region is simply connected
(see Examples \ref{ex:cubictwos} and 
\ref{ex:dumbbell}.)

\begin{figure}
\includegraphics[scale=.95]{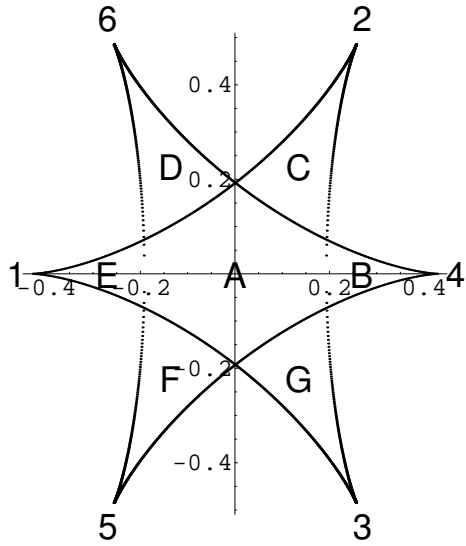}
\caption{Image of critical set for $f(z) = 
-\frac{1}{2} z^3 + \frac{3}{2} z - \overline{z}$, 
with labels}
\label{fig:cubictwos_imaglabel}
\end{figure}

\begin{figure}
\includegraphics[scale=1.50]{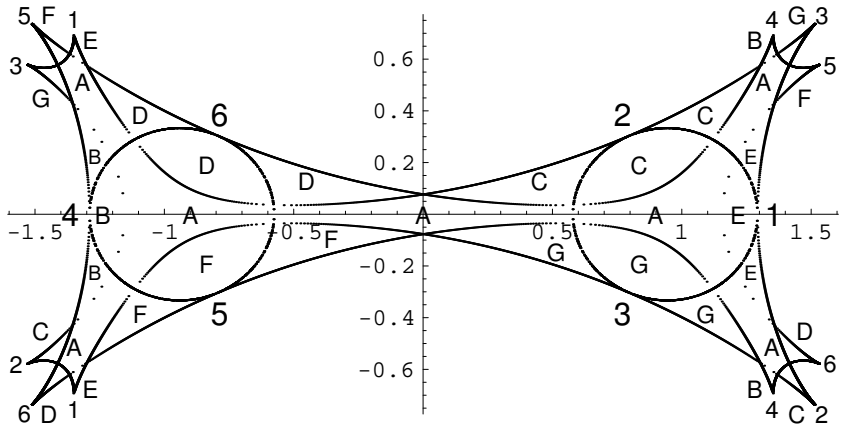}
\caption{Preimage of $f(S)$ for $f(z) =
-\frac{1}{2} z^3 + \frac{3}{2} z - \overline{z}$} 
\label{fig:cubictwos_preim}
\end{figure}

\begin{example}
\label{ex:cubictwos}
$f(z) = -\frac{1}{2} z^3 + \frac{3}{2} z - 
\overline{z}$

Clearly, $C(f, \infty) = \varnothing$;
thus the partitioning set is $f(S)$.
The critical set of $f$ is graphed in
Figure \ref{fig:cubictwos_crit};
it consists of two disjoint, closed curves.
Figure \ref{fig:cubictwos_imaglabel} shows
$f(S)$,
with the cusps numbered
and regions labeled in upper case letters.
Using Mathematica to find preimages of a few
points in each component of the partition,
we see that $Val(f, w) = 3$ on the unbounded
component,
$Val(f, w) = 5$ on the tips of the star,
and $Val(f, w) = 7$ in the ``center'' of the star.

\begin{figure}
\includegraphics[scale=.95]{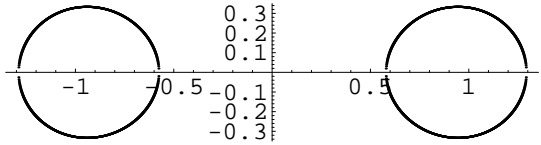}
\caption{Critical set for $f(z) = 
-\frac{1}{2} z^3 + \frac{3}{2} z - \overline{z}$}
\label{fig:cubictwos_crit}
\end{figure}

Notice that each of the image regions is simply
connected;
hence we expect $f$ to be univalent on each of
the corresponding preimage regions by 
Theorem \ref{thm:preimcover}.
This is consistent with the results
for finding the preimages of
a point in each component of 
$\mathbb{C} \backslash f(S)$ using Mathematica.

The preimages of the cusps on the critical set
are labeled in Figure \ref{fig:cubictwos_preim}
with the same numbers
used in Figure \ref{fig:cubictwos_imaglabel}.
We have labeled the preimage
components corresponding to a given bounded
component of $\mathbb{C} \backslash f(S)$ with
the corresponding letter.

As expected from Theorem \ref{thm:jointwo},
we can join components sharing a common arc
with interior off $S$ to get a larger region of
univalence.
From Figure \ref{fig:cubictwos_preim},
it is apparent that
we can continue joining such bounded components 
in the unbounded
component of $\mathbb{C} \backslash S$ 
until we reach $S$.
The combined regions in the image and in the domain
appear to be simply connected,
so univalence in the combined region should follow 
from the theorem of Ortel and Smith
(Theorem \ref{thm:OSlcluniv}.) 
\end{example} 

\begin{figure}
\includegraphics[scale=1.20]{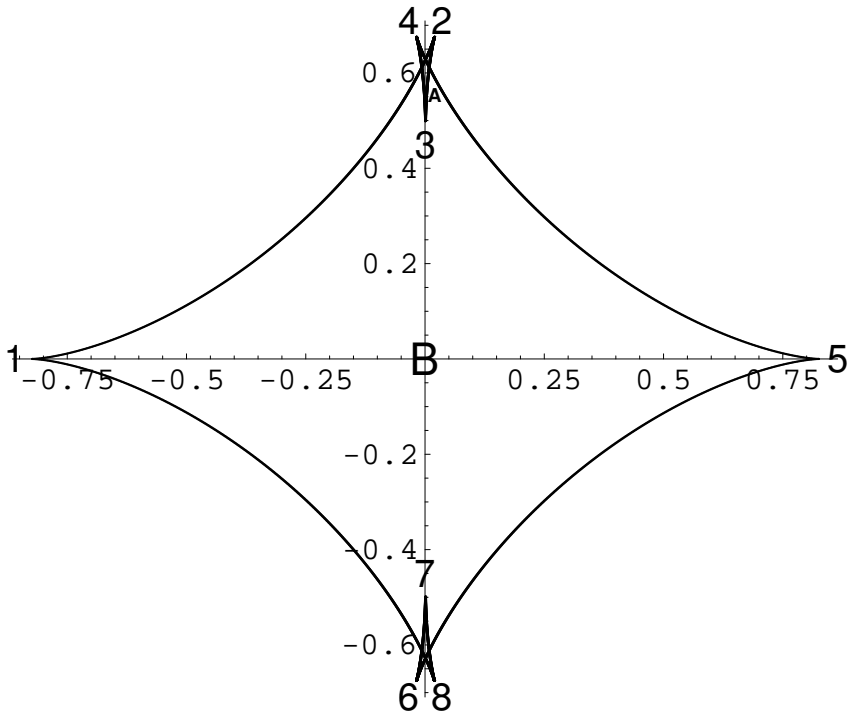}
\caption{Image of critical set for $f(z) = 
-\frac{1}{2} z^3 + \frac{9}{10}z - \overline{z}$,
with labels}
\label{fig:dumbbell_imaglabel}
\end{figure}

\begin{figure}
\includegraphics[scale=1.20]{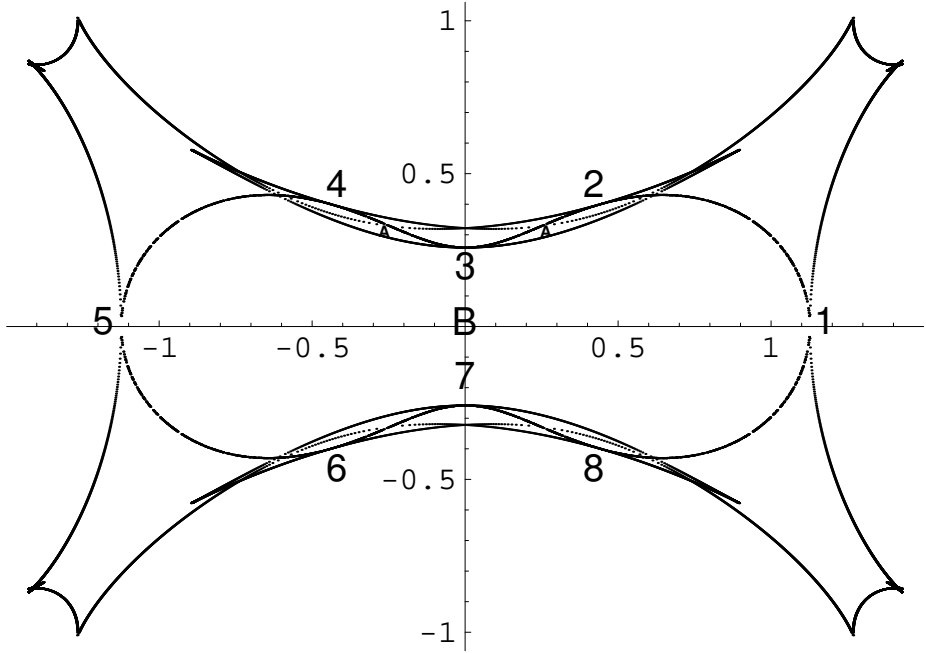}
\caption{Preimage of $f(S)$ for $f(z) =
-\frac{1}{2} z^3 + \frac{9}{10}z - \overline{z}$}
\label{fig:dumbbell_preim}
\end{figure}

\begin{example}
\label{ex:dumbbell}
$f(z) = -\frac{1}{2} z^3 + \frac{9}{10}z - 
\overline{z}$

$C(f, \infty) = \varnothing$;
thus the partitioning set is $f(S)$.
The critical set of $f$ is graphed in
Figure \ref{fig:dumbbell_crit};
it looks like a dumbbell.
Figure \ref{fig:dumbbell_imaglabel} shows
$f(S)$,
with the cusps numbered
and two regions labeled in upper case letters.
Using Mathematica to find preimages of a few
points in each component of the partition,
we see that $Val(f, w) = 3$ on the unbounded
component and increases by two each time
we cross an arc in $f(S)$ into a component
containing the tangent to the arc.
The components with maximum valence (7)
are bisected by the imaginary axis.

\begin{figure}
\includegraphics[scale=.75]{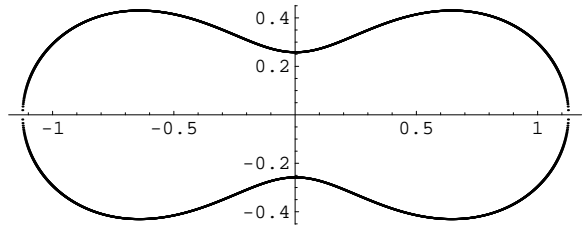}
\caption{Critical set for $f(z) = 
-\frac{1}{2} z^3 + \frac{9}{10}z - \overline{z}$}
\label{fig:dumbbell_crit}
\end{figure}

Notice that each of the image regions appears to be 
simply connected;
hence we expect $f$ to be univalent on each of
the corresponding preimage regions by 
Theorem \ref{thm:preimcover}.

In Figure \ref{fig:dumbbell_preim},
we label the regions of the
partition of the domain inside
the bounded component of $\mathbb{C} \backslash S$
which
correspond to our two labeled image regions.
Note that if we combine the three regions meeting
at the point numbered $3$
(a preimage of the cusp numbered $3$
in Figure \ref{fig:dumbbell_imaglabel}),
we run into trouble.
Let's call that point $z_0$.
Note that $z_0$ is the only shared boundary
point in $S$ for the three regions and is not
included in the combined region.
We note that $f(z_0)$ is a puncture point in
the combined region,
so the theorem of Ortel and Smith 
(Theorem \ref{thm:OSlcluniv})
cannot be applied to claim that $f$ is univalent
on the combined region in the domain.

In Figure \ref{fig:dumbbell_preim},
we only labeled the preimages in $S$ of cusps.
\end{example}
  
To make further progress with this,
we need to have a better understanding of what
the regions of the partitions look like
(for example, 
is a component of the partition of the image
ever an annulus?)
An obvious starting point is to look at 
harmonic polynomials where 
$C(f, \infty) = \varnothing$.
We make the following

\bigskip
\noindent
\textbf{Conjecture:}
\textit{
Let $f(z) = p(z) - \overline{z}$ with $p(z)$ being
a polynomial in $z$ of degree two or higher.
Let $\{R_j\}$ be a finite collection of
bounded regions in the unbounded component of 
$\mathbb{C} \backslash S$ such that each
$R_j$ is a component of $\mathbb{C} \backslash
f^{-1}(f(S))$. 
Let $R = int(\bigcup \overline{R_j})$.
Suppose that we choose the $R_j$ such that
$R$ is a connected region where
$z \in \partial R$ implies one of the following holds:
\textit{(i)} $z \in S$
or 
\textit{(ii)} every neighborhood of
$z$ contains points mapped to the unbounded component
of $\mathbb{C} \backslash f(S)$.
Then $f$ is univalent on $R$.
Moreover,
$R$ is maximal in the sense that
if we add another bounded component of
$\mathbb{C} \backslash f^{-1}(f(S))$ to $R$,
then
$f$ will not be univalent on the combined region.
}

\section {Global valence and Lyzzaik's local results}
\label{sec:lyzzaik}
The examples with $Val(f) < \infty$ where $f$
is defined in $\mathbb{C}$
in Sections \ref{sec:harm} and \ref{sec:entire} above
share a common feature:
the valence increases by an even number when one
crosses an arc of $f(S)$ into the region containing
the tangent line to that arc.
Since
$Val(f) < \infty$ in each such example, 
these are examples of 
light harmonic functions.
We shall apply some results of Lyzzaik~\cite{L:light}
concerning the behavior of a light harmonic function
defined in a simply connected domain to show that
the change in valence illustrated in these examples 
is a property of light harmonic
functions in $\mathbb{C}$.

First, we will consider the setting for Lyzzaik's
results.
Let $W$ be a simply connected domain in $\mathbb{C}$
and $f$ a harmonic function in $W$.
Then we can represent $f$ as $f = h + \overline{g}$
where $h$ and $g$ are holomorphic in $W$.
Let $\psi(z) = (h' / g')(z)$ for $z \in W$;
$\psi$ is the reciprocal of the dilatation of $f$.
Lyzzaik denotes the critical set of $f$ by $J$;
for consistency,
we will denote it by $S$.

\begin{definition}[Lyzzaik]
\label{def:lyzdef2.1}
For a light harmonic mapping $f$ in $W$
every $z \in N = \{z \in S: |\psi(z)| \ne 1\}$
is called a non-folding critical point of $f$.
Note that if $z \in N$, 
then $h'(z) = g'(z) = 0$.
\end{definition}

\begin{lemma}[Lyzzaik]
\label{lem:lyzlem2.2}
Every $z_0 \in N$ belongs to a neighborhood that
contains no other point in $S$.
$\psi$ is unimodular on $\Gamma_f = S \backslash N$ 
and $\Gamma_f$ consists of
curves which are analytic except possibly for
algebraic singularities.
\end{lemma}

\noindent
Lyzzaik parametrizes a directed Jordan subarc $\gamma$
of $\Gamma_f$ by an analytic path $z(t)$ 
for $t \in I = [0, 1]$.
He then defines a continuous, increasing function 
$\phi : I \rightarrow \mathbf{R}$
by requiring that $\psi(z(t)) =\ $exp$(i \phi(t))$.
An easy calculation shows that

\begin{equation*}
\frac{d}{dt} f(z(t)) =
2 (Re\ \omega(t)) \exp(i \phi(t) / 2)
\end{equation*}

\noindent
where

\begin{equation*}
\omega(t) = h'(z(t)) z'(t) 
\exp(-i \phi(t) / 2)
\end{equation*}

\noindent
for $z \in \gamma \subseteq S \backslash N$.

\newcounter{lyz1}
\newcounter{lyz2}
\begin{definition}[Lyzzaik]
\label{def:lyzdef2.2}
Let $f$ be a light harmonic mapping in $W$ and
let $z_0 \in \Gamma_f$.
\begin{list}{(\Alph{lyz1})}{\usecounter{lyz1}}
\item
If $\psi'(z_0) \ne 0$,
then $z_0$ is interior to a Jordan subarc $\gamma$
of $\Gamma_f$ which we assume given as above,
with $z_0 = z(t_0)$, 
$0 < t_0 < 1$.
Let
$\omega(t)$ be as above.
We call $z_0$ a critical point of $f$ of the
\begin{list}{(\alph{lyz2})}{\usecounter{lyz2}}
\item
first kind if $Re\ \omega$ changes sign at $t_0$.
\item
second kind if $z_0$ is not a critical point 
of the first
kind and if $z_0$ is a zero of $h'$,
or equivalently $g'$
(which yields $Re\ \omega(t_0)) = 0$.)
\end{list}
\item
If $\psi'(z_0) = 0$,
then we call $z_0$ a critical point of $f$ of the
third kind.
\end{list}
Let $F_j$ where $j$= 1,2, or 3,
denote the set of all critical 
points of $f$ of the
$j$th kind,
and let $F = \cup_{j=1}^3 F_j$.
We call $F$ the set of folding critical points 
of $f$.
Note that this classification is
independent of $\gamma$ and its parametrization.
\end{definition}

\begin{theorem}[Lyzzaik]
\label{thm:lyzthm2.2}
Let $f$ be a light harmonic function in $W$.
Then $F \cup N$ consists of isolated points.
\end{theorem}

\begin{remark}
\label{rem:noaccpts}
We also claim that $F \cup N$ has no point of
accumulation in $W$.
Since $f$ is light,
neither $F_3$ nor $N$ can have a point of a 
accumulation;
otherwise $S$ has nonempty interior
(in which case,
the proof of Lemma \ref{lem:nullpartifsnotempty}
implies that $f$ is not light.)
Further,
as noted by Lyzzaik,
$Re\ \omega$ can have only finitely many zeros
on any Jordan arc in $S$ when $f$ is light.
Hence,
$F_1 \cup F_2$ has no point of accumulation.
\end{remark}

We can now state one of Lyzzaik's structure
theorems for the local behavior of $f$ near
the critical set.
The notation $f_{z_0} \sim z^j, \overline{z}^k$  
is explained in Section \ref{sec:lyzintro} above.

\newcounter{lyz3}
\begin{theorem}[Lyzzaik]
\label{thm:lyzthm5.1}
Let $f = \overline{g} + h$ be a light harmonic
mapping in $W$,
$z_0 \in \Gamma_f$, and $\ell \ge 0$
the order of $z_0$ as a zero of $h'$
or equivalently $g'$.
\begin{list}{(\alph{lyz3})}{\usecounter{lyz3}}
\item
Suppose that $z_0 \in \Gamma_f \backslash
(F_1 \cup F_3)$.
Then $f$ satisfies $f_{z_0} \sim z^{\ell+1},
\overline{z}^{\ell+1}$
if $\ell$ is even (including zero),
and $f_{z_0} \sim z^{\ell + 2}, \overline{z}^\ell$
or $f_{z_0} \sim z^\ell, \overline{z}^{\ell + 2}$
if $\ell$ is odd. 
\item
Suppose that $z_0 \in F_1$.
Then $f$ satisfies
$f_{z_0} \sim z^{\ell+1},\overline{z}^{\ell+3}$
or
$f_{z_0} \sim z^{\ell+3},\overline{z}^{\ell+1}$
if $\ell$ is even,
and $f_{z_0} \sim z^{\ell + 2}, 
\overline{z}^{\ell + 2}$
if $\ell$ is odd. 
\end{list}
\end{theorem}

Lyzzaik defines a convex arc to be a directed simple 
arc that has the slope of its tangent continuously 
increasing.
An arc is locally convex at $z_0$ if $z_0$ belongs
to some open subarc which is convex.
Let $\gamma \subseteq \Gamma_f$
and $z_0 \in int\ \gamma$.
Lyzzaik examines arg\ $\frac{d}{dt} f(z(t))$ 
at $z_0 = z(t_0)$
and shows that $f(\gamma)$ is locally convex
at $z_0 \in \Gamma_f \backslash (F_1 \cup F_3)$.
Thus we can speak of the region which contains the 
tangent line
to $f(\gamma)$ at $f(z_0)$ if $z_0 \notin F$.
(Another approach is to note that if 
$\gamma \subseteq \Gamma_f$,
then $f(\gamma)$ cannot be piecewise constant;
otherwise,
a connected subset of $\gamma$ is mapped to a point
and $f$ is not light.
One can then apply a result of P. Duren
and D. Khavinson ~\cite{DK:concave}
to show that $f(\gamma)$ is concave, provided
that $f$ is univalent in a region with $\gamma$
in its boundary.
However,
this approach requires us to show that we can find
such a region of univalence.)

Let $f$ be a light harmonic function in $\mathbb{C}$.
Let $\gamma$ be a simple analytic arc in the critical
set such that $\beta = f(\gamma)$ is convex 
(and contains no points of inflection) and
such that $\gamma\ \cap\ (F \cup N)  
= \varnothing$.
An arc will be called non-degenerate if
it is a non-empty arc that does not consist of a single
point.
Here, $\gamma$ is non-degenerate.
Then, given $z_0 \in\ int\ \gamma$,
there is a neighborhood $U$ of $z_0$ such that
$\gamma$ partitions $U$ into a sense-preserving region
$U^+$
and a sense-reversing region $U^-$.
If we choose $z_0 \in int\ \gamma$,
then $\ell = 0$ in Theorem \ref{thm:lyzthm5.1}
(since $g'(z_0) \neq 0$)
and $f_{z_0} \sim z, \overline{z}$.
By Theorem \ref{thm:lyzthm5.1} and its proof,
we can choose $U$ and some neighborhood of $f(z_0)$,
say $B$,
such that $\beta$ partitions $B$
into two regions $V^+$ and $V^-$
and such that $f$ maps $U^+$ 1-1 onto $V^+$ and maps
$U^-$ 1-1 onto $V^+$.
Here, $V^+$ is chosen such that the tangent to 
$\beta$ at $f(z_0)$ lies in $V^+$.
We will use these ideas to prove:

\begin{theorem}
\label{thm:lyzval}
Let $f$ be a finite valence harmonic mapping on 
$\mathbb{C}$.
Let $\Omega_+$ and $\Omega_-$ be distinct components
of $\mathbb{C}\backslash (f(S) \cup C(f,\infty))$
that share a common non-degenerate boundary arc 
$\beta_0$.
Suppose that $int\ \beta_0$ does not lie completely
in $C(f, \infty)$.
Then there is a non-degenerate subarc $\beta$ of 
$\beta_0$
that lies outside of $C(f,\infty)$ and
such that $f^{-1}(\beta) \cap 
(F \cup N) = \varnothing$.
Choose $w_0 \in\ int\ \beta$ and choose $\Omega_+$
to be the region containing the tangent to 
$\beta$ at $w_0$.
Suppose that $w_0$ has $N_0 \ge 0$ distinct preimages
off $S$ and $N_1 > 0$ distinct preimages in $S$.
Then $Val(f, \Omega_-) = N_0$ and 
$Val(f, \Omega_+) = N_0 + 2 N_1$.
\end{theorem}
\begin{proof}
Since $f$ has finite valence, 
$f$ is light.
We wish to choose a non-degenerate, simple analytic arc
$\beta \subseteq \beta_0$ such
that $\beta \cap C(f,\infty) = \varnothing$, 
with $\beta$ in the boundary between
two regions of constant valence,
and such that $\beta \cap f(F \cup N) = \varnothing$.
The two regions are non-empty 
by Corollary \ref{cor:goodpart},
since $Val(f) < \infty$. 
By our choice of $\beta$,
a structure theorem of Lyzzaik
for light harmonic functions 
(Theorem \ref{thm:lyzthm5.1})
should help us count the distinct
preimages of points in some sufficiently small
neighborhood $B$ of $w_0$.

\begin{enumerate}

\item Existence of $\beta$:
Since $int\ \beta_0$ does not lie completely in 
$C(f,\infty)$,
we can find a bounded, closed, non-degenerate subarc 
$\tilde{\beta}$
such that $\exists\ w \in int\ \tilde{\beta}$
with $w \notin C(f,\infty)$.
Since $C(f, \infty)$ is closed,
there exists a neighborhood of $w$ 
that is disjoint from $C(f, \infty)$.
Take $\beta_1$ to be a closed,
non-degenerate subarc of $\tilde{\beta}$
contained in this neighborhood
and to have $w$ in its interior.

If $K$ is an arbitrary compact set,
then $(F \cup N) \cap K$ must be finite
by Remark \ref{rem:noaccpts}.
Since $\beta_1 \cap C(f, \infty) = \varnothing$,
the closed set $f^{-1}(\beta_1)$ is also bounded.
Thus,
$f^{-1}(\beta_1)$ is compact and
$\beta_1$ contains at most a finite number
of points with preimages in $F \cup N$.
So, we may find a bounded subarc $\beta$ 
of $\beta_0$ that neither intersects $C(f,\infty)$
nor contains points in $f(F \cup N)$.
Moreover, 
we may assume that this arc is convex and
without any points of inflection.
Since $f$ has finite valence, $f(S) \cup C(f,\infty)$
has empty interior.
By construction, 
$\beta$ will be the local 
boundary between $\Omega_+$ and $\Omega_-$.

\item Preimages of $w_0$:
Choose $w_0 \in\ int\ \beta$.
Let $f^{-1}(w_0) \cap S = \{z_1, z_2, ... , z_{N_1}\}$.
In a sufficiently small neighborhood of $z_j$,
$S$ consists of a single non-degenerate arc
since $z_j \notin N \cup F_3$ by (1).
A non-degenerate subarc of this arc is mapped to a
subarc of $\beta$.
This follows from our assumption that
$\beta_0$ locally separates components of the partition;
\textit{i.e.},
if two distinct curves in $f(S)$ cross at $w_0$,
then $\beta_0$ doesn't separate 
$\Omega_+$ and $\Omega_-$ in a sufficiently
small neighborhood of $w_0$.
Since $w_0 \in f(S)$, we have $N_1 > 0$.
Let $f^{-1}(w_0) \cap (\mathbb{C}\backslash S) 
= \{\zeta_1, \zeta_2, ... , \zeta_{N_0}\}$,
where $N_0 \ge 0$.

\item Construction of neighborhood $B$ of $w_0$:
Consider a preimage $\zeta_j$ of $w_0$ 
off of $S$.
By the inverse function theorem,
there exist $Q_j$, 
an open neighborhood of $\zeta_j$,
and $V_j$,
an open neighborhood of $w_0$,
such that $f: Q_j \rightarrow V_j$ is 1-1, onto.
Without loss of generality, 
we may assume $Q_j \subseteq
\mathbb{C} \backslash f^{-1}(f(S) \cup C(f,\infty))$. 
Now consider a preimage $z_j$ in $S$.
By construction,
$z_j \notin N \cup F_3$.
Let $\gamma$ be a simple arc in $S$ such that
$z_j \in\ int\ \gamma$
and $f(\gamma) \subseteq \beta$. 
By our choice of $\beta$,
$\gamma \cap (F \cup N) = \varnothing$.
Hence,
$g'(z_j) \neq 0$.
We may apply our interpretation above of 
Theorem \ref{thm:lyzthm5.1} 
to find an open neighborhood
$U_j$ of $z_j$ 
and an open neighborhood $W_j$ of $w_0$
such that
\begin{itemize}
\item
$\gamma$ partitions $U_j$ into a sense-preserving
region ${U_j}^+$ and a sense-reversing region
${U_j}^-$.
\item
${U_j}^+$ and ${U_j}^-$ are each mapped
1-1 onto $W_j \cap\ \Omega_+$.
\end{itemize}

Let $B_0 = \cap_{j=1}^{N_0} V_j$ and let
$B_1 = \cap_{j=1}^{N_1} W_j$.
Clearly, $B_1 \neq \varnothing$.
If $N_0 > 0$, $B_0 \neq\ \varnothing$, 
so let $B = B_0 \cap B_1$.
Otherwise, $B_0$ is empty, so let $B = B_1$.
By construction,
$B$ is open.
Finally, let 
$T = f^{-1}(B) \cap ((\cup_{j=1}^{N_0} Q_j)
\cup (\cup_{j=1}^{N_1} U_j))$.
Then every preimage of $w_0$ is in $int\ T$.
We may choose the $Q_j$ and the $U_j$
such that these sets are all pairwise disjoint.

\item Valence in $\Omega_+$:
By construction, each $w \in B \cap \Omega_+$ has 
exactly $N_0 + 2 N_1$ distinct preimages in $T$
and thus has at least
$N_0 + 2 N_1$ distinct preimages in $\mathbb{C}$.
We claim that $\exists\ w_+ \in B \cap \Omega_+$ 
such that $Val(f, w_+) = N_0 + 2 N_1$.
Suppose not.
Then $\exists\ \{w_n\} \subset B \cap \Omega_+$ such
that $w_n \rightarrow w_0$ and such that each
$w_n$ has at least $N_0 + 2 N_1 + 1$ 
distinct preimages in $\mathbb{C}$.
In particular, each $w_n$ has a preimage 
$\xi_n \notin T$.
Since $w_n \rightarrow w_0$ and 
$w_0 \notin C(f,\infty)$,
$\{\xi_n\}$ is bounded, 
so has a convergent subsequence, which we will also
denote $\{\xi_n\}$.
Let $\xi_n \rightarrow \xi_0$.
Since $f(\xi_0) = w_0$,
$\xi_0 \in int\ T$.
Thus, $\xi_n \in T$
for $n$ sufficiently large, a contradiction.
Thus, $\exists\ w_+ \in B \cap \Omega_+$ such that 
$Val(f, w_+) = N_0 + 2 N_1$.
Since,
by Theorem \ref{thm:harmpart}, 
the valence of $f$ is constant on the
region $\Omega_+$, 
$Val(f, \Omega_+) = Val(f, w_+) = N_0 + 2 N_1$.

\item Valence in $\Omega_-$:
By construction, each $w \in B \cap \Omega_-$ has 
exactly $N_0$ distinct preimages in $T$
(one preimage in each $Q_j$) 
and thus has at least
$N_0$ distinct preimages in $\mathbb{C}$.
We claim that $\exists\ w_- \in B \cap \Omega_-$ 
such that $Val(f, w_-) = N_0$.
Suppose not.
Then $\exists\ \{w_n\} \subset B \cap \Omega_-$ such
that $w_n \rightarrow w_0$ and such that each
$w_n$ has at least $N_0 + 1$ preimages
in $\mathbb{C}$.
In particular, each $w_n$ has a preimage 
$\xi_n \notin f^{-1}(B) \cap (\cup_{j=1}^{N_0} Q_j)$.
Since $w_0 \notin C(f,\infty)$,
$\{\xi_n\}$ is bounded, 
so has a convergent subsequence, which we will also
denote $\{\xi_n\}$.
Let $\xi_n \rightarrow \xi_0$.
Since $f(\xi_0) = w_0$,
$\xi_0 \in int\ T$.
Thus, for $n$ sufficiently large,
$\xi_n \in T$.
Since $f$ maps $\cup_{j=1}^{N_1} U_j$ into 
$\beta \cup \Omega_+$,
$\xi_n \in f^{-1}(B) \cap (\cup_{j=1}^{N_0} Q_j)$
for sufficiently large $n$, a contradiction.
Thus, $\exists\ w_- \in B \cap \Omega_-$ such that 
$Val(f, w_-) = N_0$.
Since the valence of $f$ is constant on the
region $\Omega_-$
(Theorem \ref{thm:harmpart} again), 
$Val(f, \Omega_-) = Val(f, w_-) = N_0$.
\end{enumerate}
\end{proof}

\begin{remark}
The proof above can also be extended to the
case where $f$ is a finite-valence,
light harmonic function 
in a simply connected, open set.
In that case,
$C(f)$ replaces $C(f, \infty)$ and
$S$ is defined only for points in our open set.
\end{remark}

\begin{example}
Consider Example \ref{ex:isolcrit} above.
We can select a few points in each region of
the partition and ask Mathematica to find the
preimages in $B(0,1)$.
$Val(f, w) = 0$ in the unbounded
component of $\mathbb{C} \backslash f(S)$.
In each ``point of the star'' region,
$Val(f, w) = 1$.
$Val(f, w) = 2$ in the ``center of the star'' region
(recall that the origin is in $f(S)$ and is not in
this region.)
The origin has one preimage in $B(0,1)$.
This is not a counterexample to this extended version 
of Theorem \ref{thm:lyzval},
because the ``star'' is $C(f)$,
not $f(S)$.
Effectively,
when we restrict $f$ to our region,
there are no points with $|z| > 1$ for $f$
to ``fold'' over $|z| = 1$.
\end{example}

\bibliographystyle{amsplain}

\begin{thebibliography} {NAM99}
\bibitem[AL 88]{AL:local}
Y. Abu-Muhanna and A. Lyzzaik,
\textit{A geometric criterion for decomposition
and multivalence},
Math. Proc. Cambridge Phil. Soc.
\textbf{103} (1988),
487--495.
MR 89e:30010.
\bibitem[Bal 91]{Balk:book}
Mark Benevich Balk,
\textit{Polyanalytic Functions},
Mathematical research, volume 63,
Akademie Verlag GmbH (1991).
MR 93k:30076.
\bibitem[BC 55]{BC:brelotchoquetlem} 
Marcel Brelot and Gustave Choquet,
\textit{Polyn\^{o}mes harmoniques et
poly\-har\-mo\-niques}, 
Second colloque sur les \'{e}quations aux
d\'{e}riv\'{e}es partielles, Bruxelles, 
1954,
pp. 45--66.  Georges Thone, Li\`{e}ge;
Masson $\&$ Cie (1955).
MR 16,1108e.
\bibitem[BHN 99]{BHN:dilat} 
Daoud Bshouty, Walter Hengartner, and 
M. Naghibi-Beidokhti, 
\textit{p-valent harmonic mappings with finite 
Blaschke dilatations},
XII-th Conference on
Analytic Functions (Lublin, 1998), 
Ann. Univ. Mariae Curie-Sklodowska Sect. A, 
\textbf{53} (1999), 
9--26.
MR 2001j:30016.
\bibitem[BHS 95]{BHS:maxval} 
Daoud Bshouty, Walter Hengartner, and
Tiferet Suez,
\textit{The exact bound on the number of zeros
of harmonic polynomials}, J. Anal. Math.
\textbf{67} (1995), 207--218.
MR 97f:30025. 
\bibitem[CL 66]{CL:cluster}
E. F. Collingwood and A. J. Lohwater, 
\textit{The Theory of
Cluster Sets}, 
Cambridge University Press (1966).
MR 38\#325.
\bibitem[DK 97]{DK:concave} 
Peter Duren and Dmitry Khavinson,
\textit{Boundary correspondence and dilatation of
harmonic mappings},
Complex Variables Theory Appl.,
\textbf{33} (1997), 105--111.
MR 98m:30039.
\bibitem[KS 03]{KS:3n-2} 
Dmitry Khavinson and 
Grzegorz \'{S}wi\c{a}tek,
\textit{On the number of zeros of certain harmonic
polynomials},
Proc. Amer. Math. Soc.
\textbf{131} (2003),
409--414.
\bibitem[Lew 36]{Lew:jacobian}
Hans Lewy,
\textit{On the non-vanishing of the Jacobian in
certain one-to-one mappings},
Bull. Amer. Math. Soc.
\textbf{42} (1936),
689--692.
\bibitem[Lyz 92]{L:light} 
Abdallah Lyzzaik, 
\textit{Local properties of light
harmonic mappings}, 
Canad. J. Math. 
\textbf{44} (1992),
135--153.
MR 93e:30048.
\bibitem[Mun 75]{Mun:top}
James R. Munkres, 
\textit{Topology: A First Course},
Prentice-Hall, Inc. (1975).
MR 57\#4063.
\bibitem[Neu 03]{gcn:thesis}
Genevra Chasanov Neumann,
Valence of harmonic functions,
Ph.D. dissertation,
University of California, Berkeley. 2003.
\bibitem[OS 86]{OrtelSmith:univalence}
M. Ortel and W. Smith,
\textit{A covering theorem for continuous
locally univalent maps of the plane},
Bull. London Math. Soc.
\textbf{18} (1986),
359--363.
MR 88b:30013.
\bibitem[Smi 71]{S:primer}
Kennan T. Smith,
\textit{Primer of Modern Analysis},
Bogden \& Quigley, Inc.
(1971).
MR 84m:26002.
\bibitem[Sto 56]{Sto:lecons}
S. St\"{o}ilow,
\textit{Le\c{c}ons sur les principes topologiques de
la the\'{o}rie des fonctions analytiques},
deuxi\`{e}me edition, Gauthier-Villars (1956).
MR 18,568b.
\bibitem[ST 00]{ST:partition}
T. J. Suffridge and J. W. Thompson,
\textit{Local behavior of harmonic mappings},
Complex Variables Theory Appl.,
\textbf{41} (2000), 63--80.
MR 2001a:30019.
\bibitem[Wil 94]{W:thesis} 
Alan Stephen Wilmshurst,
\textit{Complex harmonic mappings and the valence
of harmonic polynomials}, D.Phil. thesis,
University of York, England. 1994.
\bibitem[Wil 98]{W:paper} 
A. S. Wilmshurst, 
\textit{The valence of harmonic polynomials},
Proc. Amer. Math. Soc.
\textbf{126} (1998),
2077--2081.
MR 98h:30029.
\end{thebibliography}

\end{document}